\documentclass[11pt]{article}
\usepackage{amsmath,amssymb}
\usepackage{graphicx}
\usepackage{subfigure}
\usepackage{epsfig}
\usepackage{hyperref}

\def\sgn{{\hbox{sgn}}}
\def\Re{{\hbox{Re}}}
\def\tr{{\hbox{\rm Tr}}}
 
\def\supp{{\hbox{\rm supp}}}
\def\C{{\hbox{\bf C}}}
\def\E{{\hbox{\bf E}}}
\def\P{{\hbox{\bf P}}}
\def\Z{{\hbox{\bf Z}}}

\newenvironment{proof}{\noindent {\bf Proof} }{\endprf\par}
\def \endprf{\hfill {\vrule height6pt width6pt depth0pt}\medskip}
\def\emph#1{{\it #1}} \def\textbf#1{{\bf #1}}

\def\RR{\mathbb{R}}

\newcommand{\<}{\langle}
\renewcommand{\>}{\rangle}

\newcommand{\goto}{\rightarrow}

\newtheorem{theorem}{Theorem}[section]
\newtheorem{lemma}[theorem]{Lemma}
\newtheorem{corollary}[theorem]{Corollary}
\newtheorem{proposition}[theorem]{Proposition}
\newtheorem{definition}[theorem]{Definition}

\numberwithin{equation}{section}

\def\ZZ{\mathbb{Z}}

\def\supp{\text{supp}}
\def\cF{{\cal F}}

\begin{document}

\title{Near Optimal Signal Recovery From Random Projections:\\
Universal Encoding Strategies?}

\author{Emmanuel Candes$^{\dagger}$ and Terence Tao$^{\sharp}$\\
  \vspace{-.3cm}\\
  $\dagger$ Applied and Computational Mathematics, Caltech, Pasadena, CA 91125\\
  \vspace{-.3cm}\\
  $\sharp$ Department of Mathematics, University of California, Los
  Angeles, CA 90095
} 

\date{October 2004; Revised March 2006} 

\maketitle

\begin{abstract}  
  Suppose we are given a vector $f$ in a class ${\cal F} \subset
  \RR^N$, e.g. a class of digital signals or digital images. How many
  linear measurements do we need to make about $f$ to be able to
  recover $f$ to within precision $\epsilon$ in the Euclidean
  ($\ell_2$) metric?
  
  This paper shows that if the objects of interest are sparse in a
  fixed basis or compressible, then it is possible to reconstruct $f$
  to within very high accuracy from a small number of random
  measurements by solving a simple linear program.  More precisely,
  suppose that the $n$th largest entry of the vector $|f|$ (or of its
  coefficients in a fixed basis) obeys $|f|_{(n)} \le R \cdot
  n^{-1/p}$, where $R > 0$ and $p > 0$. Suppose that we take
  measurements $y_k = \<f, X_k\>$, $k = 1, \ldots, K$, where the $X_k$
  are $N$-dimensional Gaussian vectors with independent standard
  normal entries. Then for each $f$ obeying the decay estimate above
  for some $0 < p < 1$ and with overwhelming probability, our
  reconstruction $f^\sharp$, defined as the solution to the
  constraints $y_k = \<f^\sharp, X_k \>$ with minimal $\ell_1$ norm,
  obeys
\[
\|f - f^\sharp\|_{\ell_2} \le C_p \cdot R \cdot (K/\log N)^{-r}, \quad
r = 1/p - 1/2.
\]
There is a sense in which this result is optimal; it is generally
impossible to obtain a higher accuracy from any set of $K$
measurements whatsoever. The methodology extends to various other
random measurement ensembles; for example, we show that similar
results hold if one observes few randomly sampled Fourier coefficients
of $f$.  In fact, the results are quite general and require only two
hypotheses on the measurement ensemble which are detailed.

\end{abstract}


{\bf Acknowledgments.} E.~C. is partially supported by National
Science Foundation grants DMS 01-40698 (FRG) and ACI-0204932 (ITR),
and by an Alfred P.  Sloan Fellowship. T.~T. is supported in part by a
grant from the Packard Foundation.  E.~C.~would like to thank Justin
Romberg and Houman Ohwadi for conversations related to this project,
and Noureddine El Karoui for sending us an early draft of his
Ph.~D.~thesis.  E.~C.~would also like to thank Joel Tropp for sharing
some notes.  T. T.~thanks Wilhelm Schlag and Gerd Mockenhaupt for
pointing out the relevance of Bourgain's work \cite{bourgain-lambda},
\cite{bourgain-halasz}. Finally, we would like to thank the anonymous
referees for providing useful references (see Section
\ref{sec:precedents}) and for their role in improving the original
manuscript.

\section{Introduction and Overview of the Main Results}
\label{sec:introduction}

This paper considers the fundamental problem of recovering a finite
signal $f \in \RR^N$ from a limited set of measurements. Specifically,
given a class of signals $\cF \subset \RR^N$, one is interested in the
minimum number of linear measurements one has to make to be able to
reconstruct objects from $\cF$ to within a fixed accuracy $\epsilon$,
say, in the usual Euclidean $\ell_2$-distance. In other words, how can
one specify $K = K(\epsilon)$ linear functionals  
\begin{equation}
  \label{eq:linear-functionals}
  y_k = \<f, \psi_k\>, \quad k \in \Omega,
\end{equation}
where $(\psi_k)_{k \in \Omega}$ is a set of vectors with cardinality
$|\Omega| = K$, so that it is possible to reconstruct an object
$f^\sharp$ from the data $(y_k)_{k \in \Omega}$ obeying
\begin{equation}
  \label{eq:error}
  \|f - f^\sharp\|_{\ell_2} \le \epsilon, 
\end{equation}
for each element $f$ taken from $\cF$? The primary goal is of course,
to find appropriate functionals $(\psi_k)_{k \in \Omega}$ so that the
required number $K$ of measurements is as small as possible.  In
addition, we are also interested in concrete and practical recovery
algorithms.

The new results in this paper will address this type of question for
signals $f$ whose coefficients with respect to a fixed reference basis
obey a power-law type decay condition, and for random measurements
$(y_k)_{k \in \Omega}$ sampled from a specified ensemble.  However,
before we discuss these results, we first recall some earlier results
concerning signals of small support.  (See also Sections \ref{l1-sec} and 
\ref{connection-sec} for a more extensive discussion of related results.)

\subsection{Exact Reconstruction of Sparse Signals}
\label{sec:exact}

In a previous article \cite{CRT}, the authors together with J.~Romberg
studied the recovery of {\em sparse} signals from limited
measurements; i.e. of signals which have relatively few nonzero terms
or whose coefficients in some fixed basis have relatively few nonzero
entries. This paper discussed some surprising phenomena, and we now
review a special instance of those.  In order to do so, we first need
to introduce the discrete Fourier transform which is given by the
usual formula\footnote{Strictly speaking, the Fourier transform is
  associated to an orthonormal basis in $\C^N$ rather than $\RR^N$.
  However all of our analysis here extends easily to complex signals
  instead of real signals (except for some negligible changes in the
  absolute constants $C$).  For ease of exposition we shall focus
  primarily on real-valued signals $f \in \RR^N$, except when referring
  to the Fourier basis.}
\begin{equation}
  \label{eq:dft}
  \hat f(k) := \frac{1}{\sqrt{N}} \sum_{t \in \ZZ_N} f(t)  
e^{-i 2\pi k t/N},
\end{equation}
where the frequency index $k$ ranges over the set $\ZZ_N := \{0, 1,
\ldots, N-1\}$.  

Suppose then that we wish to recover a signal $f \in \RR^N$ made out of
$|T|$ spikes, where the set $T$ denotes the support of the signal
\[
T := \{t, \, f(t) \neq 0\}. 
\]
We do not know where the spikes are located nor do we know their
amplitudes.  However, we are given information about $f$ in the form
of `only' $K$ randomly sampled Fourier coefficients $F_\Omega f :=
(\hat f(k))_{k \in \Omega}$ where $\Omega$ is a random set of $K$
frequencies sampled uniformly at random.  In \cite{CRT}, it was shown
that $f$ could be reconstructed {\em exactly} from these data provided
that the expected number of frequency samples obeyed the lower bound
\begin{equation}
  \label{eq:T/Omega}
  |T| \le \alpha \cdot (K/\log N)
\end{equation}
for all sufficiently small $\alpha > 0$ (i.e. $\alpha \leq \alpha_0$
for some small absolute constant $\alpha_0$).  To recover $f$ from
$F_\Omega f$, we simply minimize the $\ell_1$-norm of the
reconstructed signal
\begin{equation}
  \label{eq:l1}
  \min \|g\|_{\ell_1} := \sum_{t \in \ZZ_N} |g(t)|, 
\end{equation}
subject to the constraints 
\begin{equation*}
  \hat g(k) = \hat f(k), \qquad \forall k \in \Omega. 
\end{equation*}
Moreover, the
probability that exact recovery occurs exceeds $1 - O(N^{-\rho
  /\alpha})$; $\rho > 0$ is here a universal constant and it is worth
noting that the aforementioned reference gave explicit values for this
constant.  The implied constant in the $O()$ notation is allowed to depend on $\alpha$, but
is independent of $N$.  In short, exact recovery may be achieved by solving a simple
convex optimization problem---in fact, a linear program for
real-valued signals--- which is a result of practical significance. 

In a following paper, Cand\`es and Romberg \cite{CR} extended these
results and showed that exact reconstruction phenomena hold for other
synthesis/measurement pairs.  For clarity of presentation, it will be
convenient to introduce some notations that we will use throughout the
remainder of the paper. We let $F_\Omega$ denote the $|\Omega|$ by $N$
matrix which specifies the set of those $|\Omega|$ linear functionals
which describe the measurement process so that the available
information $y$ about $f$ is of the form
\[
y = F_\Omega f.
\]
For instance, in our previous example, $F_\Omega$ is the $|\Omega|$ by
$N$ partial Fourier matrix whose rows are the sampled sinusoids
\[
F_\Omega(k,t) := \frac{1}{\sqrt{N}} \, e^{-i2\pi k t/N}, \quad k \in
\Omega, \, t \in \ZZ_N.
\]
More generally, suppose that one is given an orthonormal basis $\Psi$
\[
\Psi = (\psi_k(t))_{0 \le t,k < N},
\] 
and that one has available partial information about $f$ in the sense
that we have knowledge about a randomly selected set $\Omega \subset
\{0,\ldots,N-1\}$ of coefficients in basis $\Psi$. For convenience,
define $\Psi$ to be the $N$ by $N$ synthesis matrix with entries
$\Psi(t,k) := \psi_k(t)$.  Then $F_\Omega$ is now obtained from
$\Psi^*$ by extracting the $|\Omega|$ rows with indices $k$ obeying $k
\in \Omega$.  Suppose as before that there is another (fixed)
orthonormal basis $\Phi$ in which the coefficients $\theta(f) =
(\theta_t(f))_{1 \leq t\leq N}$ of $f$ in this basis, defined by
\[
\theta_t (f) := \<f, \phi_t\>, 
\]
are sparse in the sense that only few of the entries of $\theta(f)$
are nonzero. Then it was shown in \cite{CR} that with overwhelming
probability, $f$ is the solution to
\begin{equation}
   \label{eq:nP1}
    \min_g \|\theta(g)\|_{\ell_1} \qquad \text{subject to} 
 \qquad F_\Omega \, g = F_\Omega f.
 \end{equation}
 That is, exact reconstruction still occurs; the relationship here
 between the number of nonzero terms in the basis $\Phi$ and the
 number of observed coefficients $|\Omega|$ depends upon the {\em
   incoherence} between the two bases.  The more incoherent, the fewer
 coefficients needed; in the other direction, in the maximally
 coherent case, e.g. $\Psi = \Phi$, one in fact needs to sample
 essentially all of the coefficients in order to ensure exact
 reconstruction (the same holds true if $\Phi$ and $\Psi$ share only
 one element with nonzero inner product with $f$).
 
 A special instance of these results concerns the case where the set
 of measurements is generated completely at random; that is, we sample
 a random orthonormal basis of $\RR^N$ and observe only the first $K$
 coefficients in that basis (note that there is no advantage in
 randomizing $\Omega$ as in Section \ref{sec:exact} since the basis is
 already completely random). As before, we let $F_\Omega$ be the
 submatrix enumerating those sampled vectors and solve \eqref{eq:nP1}.
 Then a consequence of the methodology developed in this paper is that
 exact reconstruction occurs with probability at least $1 - O(N^{-\rho
   /\alpha})$ (for a different value of $\rho$) provided that
 \begin{equation}
   \label{eq:size-constraints}
   \|\theta(f)\|_{\ell_0}  \le \alpha \cdot (K/\log N), 
 \end{equation}
 where $\alpha > 0$ is sufficiently small, and
 the $\ell_0$-norm is of course the size of the support of the
 vector $\theta$
\[
\|\theta\|_{\ell_0} := |\{ t : \theta_t \neq 0\}|,
\]
see \cite{Equivl0l1} for sharper results.  In summary, $\ell_1$
  seems to recover sparse unknown signals in a variety of different
  situations. The number of measurements simply needs to exceed the
  number of unknown nonzero coefficients by a proper amount.

Observe that a nice feature of the random basis discussed above is its
statistical invariance by rotation. Let $\Phi$ be any basis so that
$\theta(f)$ are the coefficients of $f$ in that basis: $\theta(f) :=
\Phi^* f$. The constraints in \eqref{eq:nP1} impose
\[
F_\Omega \Phi \,\, \theta(g)  = F_\Omega \Phi \,\, \theta(f)
\]
and since the distribution of $F_\Omega \Phi$ is that of $F_\Omega$,
the choice of the basis $\Phi$ is actually irrelevant. Exact
reconstruction occurs (with overwhelming probability) when the signal
is sparse in any fixed basis; of course, the recovery algorithm
requires knowledge of this basis.

\subsection{Power laws}

In general, signals of practical interest may not be supported in
space or in a transform domain on a set of relatively small size.
Instead, the coefficients of elements taken from a signal class decay
rapidly, typically like a power law \cite{DevoreActa,Donoho-IEEE}.  We
now give two examples leaving mathematical rigor aside in the hope of
being more concise.
\begin{itemize}
\item {\em Smooth signals.} It is well-known that if a continuous-time
  object has $s$ bounded derivatives, then the $n$th largest entry of
  the wavelet or Fourier coefficient sequence is of size about
  $1/n^{s+1/2}$ in one dimension and more generally, $1/n^{s/d+1/2}$
  in $d$ dimensions \cite{Donoho-IEEE}.  Hence, the decay of
  Fourier or wavelet coefficients of smooth signals exhibits a power
  law.
  
\item {\em Signals with bounded variations.} A popular model for
  signal/analysis is the space of objects with bounded variations. At
  the level of the continuum, the total-variation norm of an object is
  approximately the $\ell_1$ norm of its gradient. In addition, there
  are obvious discrete analogs for finite signals where the gradient
  is replaced by finite differences.  Now a norm which is almost
  equivalent to the total-variation norm is the weak-$\ell_1$ norm in
  the wavelet domain; that is, the reordered wavelet coefficients of a
  compactly supported object $f$ approximately decay like $1/n$
  \cite{CohenPetrushev}. At the discrete level, $\|f\|_{BV}$ essentially
  behaves like the $\ell_1$-norm of the Haar wavelet coefficients up
  to a multiplicative factor of at most $\log N$. Moreover, it is
  interesting to note that studies show that the empirical wavelet
  coefficients of photographs of natural scenes actually exhibit the
  $1/n$-decay \cite{MallatWaveTour}.
\end{itemize}
In fact, finding representations with rapidly decaying coefficients is
a very active area of research known as Computational Harmonic
Analysis and there are of course many other such examples. For
instance, certain classes of oscillatory signals have rapidly decaying
Gabor coefficients \cite{Feichtinger}, certain types of images with
discontinuities along edges have rapidly decaying curvelet
coefficients \cite{CurveEdge} and so on.


Whereas \cite{CRT} considered signals $f$ of small support, we now
consider objects whose coefficients in some basis decay like a
power-law. We fix an orthonormal basis $\Phi = (\phi_t)_{1 \leq t \leq N}$
(which we call the \emph{reference basis}),
and rearrange the entries $\theta_t(f) := \langle f, \phi_t \rangle$
of the coefficient vector $\theta(f)$ in decreasing order of magnitude
$|\theta|_{(1)} \ge |\theta|_{(2)} \ge \ldots \ge |\theta|_{(N)}$.  We
say that $\theta(f)$ belongs to the weak-$\ell_p$ ball or radius $R$
(and we will sometimes write $f \in w\ell_p(R)$) for some $0 < p < \infty$ and
$C > 0$ if for each $1
\le n \le N$,
\begin{equation}\label{weak}
|\theta|_{(n)} \le R \cdot n^{-1/p}.
\end{equation}
In other words, $p$ controls the speed of the decay: the smaller $p$,
the faster the decay.  The condition \eqref{weak} is also equivalent
to the estimate
$$
|\{t \in \ZZ_N : |\theta_t(f)| \geq \lambda \}| \leq
\frac{R^p}{\lambda^p} $$
holding for all $\lambda > 0$.  We shall
focus primarily on the case $0 < p < 1$.

It is well-known that the decay rate of the coefficients of $f$ is
linked to the `compressibility' of $f$, compare the widespread use of
transform coders in the area of lossy signal or image compression.
Suppose for instance that all the coefficients $(\theta_t(f))_{1 \le n
  \le N}$ are known and consider the partial reconstruction
$\theta_K(f)$ (where $1 \leq K \leq N$ is fixed) obtained by keeping
the $K$ largest entries of the vector $\theta(f)$ (and setting the
others to zero). Then it immediately follows from \eqref{weak} that
the approximation error obeys
\[
\|\theta(f) - \theta_K(f)\|_{\ell_2} \le C_p \cdot R \cdot K^{-r}, \qquad r :=
1/p - 1/2,
\]
for some constant $C_p$ which only depends on $p$.  And thus, it
follows from Parseval that the approximate signal $f_K$ obtained by
keeping the largest coefficients in the expansion of $f$ in the
reference basis $\Phi$ obeys the same estimate, namely,
\begin{equation}
  \label{eq:best-mterm}
  \|f - f_K\|_{\ell_2} \le C_p \cdot R \cdot K^{-r},
\end{equation}
where $C_p$ only depends on $p$. 

\subsection{Recovery of objects with power-law decay}

We now return to the setup we discussed earlier, where we select $K$
orthonormal vectors $\psi_1, \ldots, \psi_K$ in $\RR^N$ uniformly at
random. Since applying a fixed orthonormal transform does not change
the problem, we may just as well assume that $\Phi$ is the identity
and solve
\begin{equation}
  \label{eq:P1}
(P_1)  \qquad 
\min_{g \in \RR^N} \, \|g\|_{\ell_1} \quad  \text{subject to} 
\quad F_\Omega \, g = F_\Omega f,   
\end{equation}
where as usual, $F_\Omega f = (\langle f, \psi_k \rangle)_{k \in
  \Omega}$. In the setting where $f$ does not have small support, we
do not expect the recovery procedure \eqref{eq:P1} to recover $f$
exactly, but our first main theorem asserts that it will recover $f$
\emph{approximately}. 

{\bf Note.} From now on and for ease of exposition, we will take
$(P_1)$ as our abstract recovery procedure where it is understood that
$f$ is the sparse object of interest to be recovered; that is, $f$
could be a signal in $\RR^N$ or its coefficients in some fixed basis
$\Phi$.


\begin{theorem}
  [Optimal recovery of $w\ell_p$ from random
  measurements]\label{lp-control} Suppose that $f \in \RR^N$ obeys
  \eqref{weak} for some fixed $0 < p < 1$ or $\|f\|_{\ell_1} \le R$
  for $p = 1$, and let $\alpha > 0$ be a sufficiently small number
  (less than an absolute constant).  Assume that we are given $K$
  random measurements $F_\Omega f$ as described above.  Then with
  probability 1, the minimizer $f^\sharp$ to \eqref{eq:P1} is unique.
  Furthermore, with probability at least $1 - O(N^{-\rho/\alpha})$, we
  have the approximation
\begin{equation}
\label{approximation}
\|f - f^\sharp \|_{\ell_2} \le C_{p,\alpha} \cdot R \cdot (K/\log N)^{-r}, 
\quad r = 1/p - 1/2. 
\end{equation}
Here, $C_{p,\alpha}$ is a fixed constant depending on $p$ and $\alpha$
but not on anything else. The implicit constant in
$O(N^{-\rho/\alpha})$ is allowed to depend on $\alpha$.
\end{theorem}
The result of this theorem may seem surprising. Indeed,
\eqref{approximation} says that if one makes $O(K \log N)$ random
measurements of a signal $f$, and then reconstructs an approximate
signal from this limited set of measurements in a manner which
requires no prior knowledge or assumptions on the signal (other than
it perhaps obeys some sort of power law decay with unknown parameters)
\emph{one still obtains a reconstruction error which is equally as
  good as that one would obtain by knowing \emph{everything} about $f$
  and selecting the $K$ largest entries of the coefficient vector
  $\theta(f)$}; thus the amount of ``oversampling'' incurred by this
random measurement procedure compared to the optimal sampling for this
level of error is only a multiplicative factor of $O(\log N)$. To
avoid any ambiguity, when we say that no prior knowledge or
information is required about the signal, we mean that the
reconstruction algorithm does not depend upon unknown quantities such
as $p$ or $R$.

Below, we will argue that we cannot, in general, design a set of $K$
measurements that would allow essentially better reconstruction errors
by any method, no matter how intractable.  As we will see later,
Theorem \ref{lp-control} is a special case of Theorem
\ref{general-lp-control} below (but for the uniqueness claim which is
proved in Section \ref{sec:l2}).

\subsection{Precedents}
\label{sec:precedents}

A natural question is whether the number of random samples we
identified in Theorem \ref{general-lp-control} is, in some sense,
optimal. Or would it be possible to obtain similar accuracies with far
fewer observations?  To make things concrete, suppose we are
interested in the recovery of objects with bounded $\ell_1$-norm,
e.g. the $\ell_1$-ball
\[
{\cal B}_1 := \{f : \|f\|_{\ell_1} \le 1\}. 
\] 
Suppose we can make $K$ linear measurements about $f\in {\cal B}_1$
of the form $y = F_\Omega f$.  Then what is the best
measurement/reconstruction pair so that the error
\begin{equation}
  \label{eq:best}
  E_K({\cal B}_1) = \sup_{f \in B_1} \|f - D(y)\|_2, \quad y =  F_\Omega f, 
\end{equation}
is minimum? In \eqref{eq:best}, $D$ is the reconstruction algorithm. 

To develop an insight about the intrinsic difficulty of our
problem, consider the following geometric picture. Suppose we take $K$
measurements $F_\Omega f$; this says that $f$ belongs to 
an affine space $f_0 + S$ where $S$ is a linear subspace of
co-dimension less or equal to $K$.  Now {\em the data available for
the problem cannot distinguish any object belonging to that plane.}
Assume $f$ is known to belong to the $\ell_1$-ball ${\cal B}_1$, say, then
the data cannot distinguish between any two points in the intersection
${\cal B}_1 \cap f_0 + S$.  Therefore, any reconstruction procedure $f^*(y)$
based upon $y = F_\Omega f$ would obey
\begin{equation}
\label{eq:diameter}
\sup_{f \in \cF} \|f - f^*\| \ge \frac{\text{diam}({\cal B}_1 \cap S)}{2}. 
\end{equation}
(When we take the supremum over all $f$, we may just assume that $f$
be orthogonal to the measurements $(y = 0)$ since the diameter will of
course be maximal in that case.)  The goal is then to find $S$ such
that the above diameter be minimal.  This connects with the agenda of
approximation theory where this problem is known as finding the
Gelfand $n$-width of the class ${\cal B}_1$ \cite{PinkusBook}, as we
explain below.

The {\em Gelfand} numbers of a set $\cF$ are defined as
\begin{equation}
\label{eq:Gelfand-number}
d_K(\cF) = \inf_S \{\sup_{f \in \cF} \|P_S f\| : \text{codim}(S) < K\},
\end{equation}
where $P_S$ is, of course, the orthonormal projection on the subspace
$S$. Then it turns out that $d_K(\cF) \le E_K(\cF) \le d_K(\cF)$. Now
a seminal result of Kashin \cite{Kashin} and improved by Garnaev and
Gluskin \cite{GarnaevGluskin,Ball} shows that for the $\ell_1$ ball,
the Gelfand numbers obey
\begin{equation}
\label{eq:Kashin}
  C \, \sqrt{\frac{\log(N/K)+1}{K}} \le  d_k({\cal B}_1) \le  
C' \, \sqrt{\frac{\log(N/K)+1}{K}}, 
\end{equation}
where $C$, $C'$ are universal constants. Gelfand numbers are also
approximately known for weak-$\ell_p$ balls as well. 

Viewed differently, Kashin, Garnaev and Gluskin assert that with $K$
measurements, the minimal reconstruction error \eqref{eq:best} one can
hope for is bounded below by a constant times
$(K/\log(N/K))^{-1/2}$. In this sense, Theorem
\ref{general-lp-control} is optimal (within a multiplicative constant)
at least for $K \asymp N^\beta$, with $\beta < 1$\footnote{Note added
  in proof: since submission of this paper, we proved in
  \cite{DecodingLP} that Theorem \ref{general-lp-control} holds with
  $\log(N/K)$ instead of $\log N$ in \eqref{approximation}.}.  Kashin
also shows that if we take a random projection, $\text{diam}({\cal
  B}_1 \cap S$ is bounded above by the right-hand side of
\eqref{eq:Kashin}. We would also like to emphasize that similar types
of recovery have also been known to be possible in the literature of
theoretical computer science, at least in principle, for certain types
of random measurements \cite{noga}.  On the one hand, finding the
Chebyshev center of $\text{diam}({\cal B}_1 \cap S)$ is a convex
problem, which would yield a near-optimal reconstruction algorithm. On
the other hand, this problem is computationally intractable when $p <
1$. Further, one would need to know $p$ and the radius of the
weak-$\ell_p$ ball which is not realistic in practical applications.

The novelty here is that the information about $f$ can be retrieved
from those random coefficients by minimizing a simple linear program
\eqref{eq:P1}, and that the decoding algorithm adapts automatically to
the weak-$\ell_p$ signal class, without knowledge thereof. Minimizing
the $\ell_1$-norm gives nearly the best possible reconstruction error
simultaneously over a wide range of sparse classes of signals; no
information about $p$ and the radius $R$ are required. In addition and
as we will see next, another novelty is the general nature of the
measurement ensemble.

It should also be mentioned that when the measurement ensemble
consists of Fourier coefficients on a random arithmetic progression, a
very fast recovery algorithm that gives near-optimal results for
arbitrary $\ell_2$ data has recently been given in
\cite{GilbertStraussII}. Since the preparation of this manuscript, we
have learnt that results closely related to those in this paper have
appeared in \cite{CompressedSensing}.  We compare our results with
both these works in Section \ref{connection-sec}.

\subsection{Other Measurement Ensembles}

Underlying our results is a powerful machinery essentially relying on
properties of random matrices which gives us very precise tools
allowing to quantify how much of a signal one can reconstruct from
random measurements. In fact, Theorem \ref{lp-control} holds for other
{\em measurement ensembles.}  For simplicity, we shall consider three
types of measured data:
\begin{itemize}
\item {\em The Gaussian ensemble:} Here, we suppose that $1 \leq K
  \leq N$ and $\Omega := \{1,\ldots,K\}$ are fixed, and the entries of
  $F_\Omega$ are identically and independently sampled from a standard
  normal distribution
\[
F_\Omega(k,t) := \frac{1}{\sqrt{N}} X_{kt}, \quad  X_{kt}  \,\, i.i.d. \,\, 
N(0,1). 
\]
The Gaussian ensemble is invariant by rotation since for any fixed
orthonormal matrix $\Phi$, the distribution of $F_\Omega$ is that of
$F_\Omega \Phi$.  

\item {\em The binary ensemble:} Again we take $1 \leq K \leq N$ and
  $\Omega := \{1,\ldots,K\}$ to be fixed.  But now we suppose that the
  entries of $F_\Omega$ are identically and independently sampled from
  a symmetric Bernoulli distribution
\[
F_\Omega(k,t) := \frac{1}{\sqrt{N}} X_{kt}, \quad  X_{kt}  \,\, i.i.d. \,\, 
P(X_{kt} = \pm 1) = 1/2. 
\]

\item {\em The Fourier ensemble:} This ensemble was discussed earlier,
  and is obtained by randomly sampling rows from the orthonormal $N$
  by $N$ Fourier matrix ${\cal F}(k,t) = \exp(-i 2\pi k
  t/N)/\sqrt{N}$.  Formally, we let $0 < \tau < 1$ be a fixed
  parameter, and then let $\Omega$ be the random set defined 
by  $$\Omega = \{k : I_k = 1\},$$ where the $I_k$'s are i.i.d.~Bernoulli variables  with $\P(I_k = 1) = \tau$. We then let $R_\Omega: \ell_2(\ZZ_N) \to \ell_2(\Omega)$ be the
  restriction map $(R_\Omega g)(k) = g(k)$ for all $k \in \Omega$ (so
  that the adjoint $R^*_\Omega: \ell_2(\Omega) \to \ell_2(\ZZ_N)$ is
  the embedding obtained by extending by zero outside of $\Omega$),
  and set
\[
F_\Omega := R_\Omega {\cal F}. 
\]
In this case, the role of $K$ is played by the quantity $K :=
\E(|\Omega|) = \tau N$.  (In fact $|\Omega|$ is usually very close to
$K$; see Lemma \ref{ms-lemma}).
\end{itemize}
Just as Theorem \ref{lp-control} suggests, this paper will show that
it is possible to derive recovery rates for all three measurement
ensembles. The ability to recover a signal $f$ from partial random
measurements depends on key properties of those measurement ensembles
that we now discuss.

\subsection{Axiomatization}

We shall now unify the treatment of all these ensembles by considering
an abstract \emph{measurement matrix} $F_\Omega$, which is a random
$|\Omega| \times N$ matrix following some probability distribution
(e.g. the Gaussian, Bernoulli, or Fourier ensembles).  We also allow
the number of measurements $|\Omega|$ to be a random variable taking
values between $1$ and $N$, and set $K := \E(|\Omega|)$--- the
expected number of measurements.  For ease of exposition we shall
restrict our attention to real-valued matrices $F_\Omega$; the
modifications required to cover complex matrices such as those given
by the Fourier ensemble are simple.  We remark that we do not assume
that the rows of the matrix $F_\Omega$ form an orthogonal family.

This section introduces two key properties on $F_{\Omega}$ which---if
satisfied---will guarantee that the solution to the problem
\eqref{eq:P1} will be a good approximation to the unknown signal $f$
in the sense of Theorem \ref{lp-control}.

First, as in \cite{CRT}, our arguments rely, in part, on the
quantitative behavior of the singular values of the matrices $
F_{\Omega T}: = F_\Omega R_T^* : \ell_2(T) \goto \ell_2(\Omega) $
which are the $|\Omega|$ by $|T|$ matrices obtained by extracting
$|T|$ columns from $F_\Omega$ (corresponding to indices in a set $T$).
More precisely, we shall need to assume the following hypothesis
concerning the minimum and maximum eigenvalues of the square matrix
$F_{\Omega T}^* F_{\Omega T}: \ell_2(T) \to \ell_2(T)$.
\begin{definition}[{\bf UUP}: Uniform Uncertainty Principle]
\label{def:uup}  
We say that a measurement matrix $F_\Omega$ obeys the uniform
uncertainty principle with oversampling factor $\lambda$ if for every
sufficiently small $\alpha > 0$, the following statement is true with
probability at least\footnote{Throughout this paper, we allow implicit constants in the $O()$
  notation to depend on $\alpha$.} $1 - O(N^{-\rho/\alpha})$ for some
fixed positive constant $\rho > 0$: for \emph{all} subsets $T$ such
that
\begin{equation}
    \label{eq:TK}
    |T| \le \alpha \cdot K/\lambda,
  \end{equation}
  the matrix $F_{\Omega T}$ obeys the bounds
\begin{equation}
  \label{eq:uup}
  \frac{1}{2} \cdot \frac{K}{N} \le \lambda_{min}(F_{\Omega T}^* F_{\Omega T}) 
  \le \lambda_{max}(F_{\Omega T}^* F_{\Omega T}) 
\le \frac{3}{2} \cdot \frac{K}{N}. 
\end{equation}
Note that \eqref{eq:uup} is equivalent to the inequality
\begin{equation}\label{uup-2}
 \frac{1}{2} \frac{K}{N} \| f \|_{\ell_2}^2
\leq \| F_{\Omega} f \|_{\ell_2}^2 \leq \frac{3}{2} \frac{K}{N} \|f\|_{\ell_2}^2
\end{equation}
holding for all signals $f$ with support size less or equal to $\alpha
K/\lambda$.
\end{definition}
There is nothing special about the constants $1/2$ and $3/2$ in
\eqref{eq:uup}, which we merely selected to make the ${\bf UUP}$ as
concrete as possible. Apart from the size of certain numerical constants (in particular, implied constants 
in the $O()$ notation), nothing in our arguments depends on this special choice, and we could replace the 
pair $(1/2, 3/2)$ with a pair $(a,b)$
where $a$ and $b$ are bounded away from zero and infinity. This remark
is important to keep in mind when we will discuss the {\bf UUP} for
binary matrices.

To understand the content of \eqref{eq:uup}, suppose that $F_\Omega$
is the partial Fourier transform and suppose we have a signal $f$
supported on a set $T$ obeying $|T| \leq \alpha K/\lambda$. Then
\eqref{eq:uup} says that $\|\hat f\|_{\ell_2(\Omega)}$ is at most
$\sqrt{3K/2N} \|f\|_{\ell_2}$ with overwhelming probability.
Comparing this with Plancherel's identity $\|\hat f\|_{\ell_2(\ZZ_N)}
= \|f\|_{\ell_2}$, we see that (with overwhelming probability) a
sparse signal $f$ cannot be concentrated in frequency on $\Omega$
\emph{regardless of the exact support of $f$}, unless $K$ is
comparable to $N$.  This justifies the terminology ``Uncertainty
Principle''.  A subtle but crucial point here is that, with
overwhelming probability, we obtain the estimate \eqref{eq:uup} for
\emph{all} sets $T$ obeying \eqref{eq:TK}; this is stronger than
merely asserting that each set $T$ obeying \eqref{eq:TK} obeys
\eqref{eq:uup} separately with overwhelming probability, since in the
latter case the number of sets $T$ obeying \eqref{eq:TK} is quite
large and thus the union of all the exceptional probability events
could thus also be quite large.  This justifies the terminology
``Uniform''.  As we will see in Section \ref{sec:l2}, the uniform
uncertainty principle hypothesis is crucial to obtain estimates about
the $\ell_2$ distance between the reconstructed signal $f^\sharp$ and
the unknown signal $f$.

The {\bf UUP} is similar in spirit to several standard principles and
results regarding random projection, such as the famous
Johnson-Lindenstrauss lemma \cite{john} regarding the preservation of
distances between a finite number of points when randomly projected to
a medium-dimensional space.  There are however a number of notable
features of the {\bf UUP} that distinguish it from more standard
properties of random projections.  Firstly, there is a wide latitude
in how to select the measurement ensemble $F_\Omega$; for instance,
the entries do not have to be independent or Gaussian, and it is even
conceivable that interesting classes of completely deterministic
matrices obeying the {\bf UUP} could be constructed.  Secondly, the
estimate \eqref{eq:uup} has to hold for \emph{all} subsets $T$ of a
certain size; for various reasons in our applications, it would not be
enough to have \eqref{eq:uup} merely on an overwhelming proportion of
such sets $T$.  This makes it somewhat trickier for us to verify the
{\bf UUP}; in the Fourier case we shall be forced to use some entropy
counting methods of Bourgain.

We now introduce a second hypothesis (which appears implicitly in
\cite{CRT}, \cite{CR}) whose significance is explained below.
\begin{definition}[ERP: Exact Reconstruction Principle]
\label{def:ERP}  
We say that a measurement matrix $F_\Omega$ obeys the exact
reconstruction principle with oversampling factor $\lambda$ if for all sufficiently small
$\alpha >0$, 
each fixed subset $T$ obeying \eqref{eq:TK} and each `sign' vector
$\sigma$ defined on $T$, $|\sigma(t)| = 1$, there exists with
overwhelmingly large probability a vector $P \in \RR^N$ with the
following properties:
\begin{itemize}
\item[(i)] $P(t) = \sigma(t)$, for all  $t \in T$;
\item[(ii)] $P$ is a linear combination of the rows of $F_\Omega$
  (i.e. $P = F_\Omega^* V$ for some vector $V$ of length $|\Omega|$);
\item[(iii)] and $|P(t)| \leq \frac{1}{2}$ for all $t \in T^c :=
  \{0,\ldots,N-1\} \backslash T$.
\end{itemize}
By `overwhelmingly large', we mean that the probability be at least $1
- O(N^{-\rho/\alpha})$ for some fixed positive constant $\rho > 0$ (recall that the 
implied constant is allowed to depend on $\alpha$).
\end{definition}

Section \ref{sec:l1} will make clear that {\bf ERP} is crucial to
check that the reconstruction $f^\sharp$ is close, in the
$\ell_1$-norm, to the vector obtained by truncating $f$, keeping only
its largest entries.  Note that, in contrast to the {\bf UUP}, in {\bf
  ERP} we allow a separate exceptional event of small probability for
\emph{each} set $T$, rather than having a uniform event of high
probability that covers all $T$ at once.  There is nothing special
about the factor $1/2$ in $(iii)$; any quantity $\beta$ strictly
between 0 and 1 would suffice here.

To understand how {\bf ERP} relates to our problem, suppose that $f$
is a signal supported on a set $T$. Then using duality theory,
it was shown in \cite{CRT} that the solution to \eqref{eq:P1} is exact if and
only if there exist a $P$ with the above properties for $\sigma(t) =
\sgn(f)(t)$---hence the name.

The hypotheses {\bf UUP} and {\bf ERP} are closely related.  For
instance, one can use {\bf UUP} to prove a statement very similar to
{\bf ERP}, but in the $\ell_2$ norm rather than the $\ell_\infty$
norm; see Corollary \ref{choose} below.  One also has an implication
of the form {\bf UUP} $\implies$ {\bf ERP} for \emph{generic} signals
$f$ assuming an additional weaker hypothesis {\bf WERP}, see Section
\ref{sec:low-level}.  In \cite{CRT} and in \cite{CR}, the property
{\bf UUP} was used (together with some additional arguments) to deduce\footnote{Note added in proof: In a sequel \cite{DecodingLP} to this paper, we show that a slight strengthening of the ${\bf UUP}$ (in which the constants $\frac{1}{2}$ and $\frac{3}{2}$ are replaced by other numerical constants closer to $1$) in fact implies ${\bf ERP}$ unconditionally.}
{\bf ERP}.

We now are in position to state the main result of this paper.

\begin{theorem}\label{general-lp-control} 
  Let $F_\Omega$ be a measurement process such that {\bf UUP} and
  {\bf ERP} hold with oversampling factors $\lambda_1$ and
  $\lambda_2$ respectively.  Put $\lambda = \max(\lambda_1,
  \lambda_2)$ and assume $K \geq \lambda$. Suppose that $f$ is a
  signal in $\RR^N$ obeying \eqref{weak} for some fixed $0 < p < 1$ or
  $\|f\|_{\ell_1} \le R$ for $p = 1$, and let $r := 1/p - 1/2$.  Then
  for any sufficiently small $\alpha$,
  any minimizer $f^\sharp$ to \eqref{eq:P1} will obey
\begin{equation}
\label{approximation-general}
 \|f - f^\sharp \|_{\ell_2} \le C_{p,\alpha} \cdot R \cdot (K/\lambda)^{-r}
\end{equation}
with probability at least $1 - O(N^{-\rho/\alpha})$. The implied constant may
depend on $p$ and $\alpha$ but not on anything else.
\end{theorem}
In this paper, we will show that the Gaussian and binary ensembles
mentioned earlier obey {\bf UUP} and {\bf ERP} with $\lambda = \log
N$, while the Fourier ensemble obeys {\bf UUP} with $\lambda = (\log
N)^6$ and {\bf ERP} with $\lambda = \log N$. Hence given an object $f
\in w\ell_p(R)$, we prove that if we collect $K \geq \log N$ Gaussian
or binary measurements, then
\begin{equation}
  \label{eq:gauss-binary}
  \|f - f^\sharp \|_{\ell_2} \le O(1) \cdot R \cdot (K/\log N)^{-r}
\end{equation}
except for a set of probability at most $O(N^{-\rho/\alpha})$.  For
randomly sampled frequency data (with at least $(\log N)^6$
frequencies being sampled), the quality of the reconstruction now
reads as
\begin{equation}
  \label{eq:freq}
  \|f - f^\sharp\|_{\ell_2} \le  O(1) \cdot R \cdot (K/(\log N)^6)^{-r}.
\end{equation}

We prove this theorem in Section \ref{gen-sec}.
Observe that our earlier Theorem \ref{lp-control} follows from
\eqref{eq:gauss-binary} and is thus a special case of Theorem
\ref{general-lp-control}. Indeed, for a fixed $F_\Omega$, \eqref{eq:P1} is
equivalent to
\[
  \min_g \|\theta(g)\|_{\ell_1} \qquad \text{subject to} 
 \qquad P_\Omega g = P_\Omega f.
\] 
where $P_\Omega$ is the orthogonal projection onto the span of the
rows of $F_\Omega$. Now suppose as in the Gaussian ensemble that
$F_\Omega$ is a matrix with i.i.d. $N(0,1/N)$ entries, then $P_\Omega$
is simply the projection onto a random plane of dimension $K$ (with
probability 1) which, of course, is the setup of Theorem
\ref{lp-control}.

\subsection{About the $\ell_1$ norm}\label{l1-sec}

We would like to emphasize that the simple nonlinear reconstruction
strategy which minimizes the $\ell_1$-norm subject to consistency with
the measured observations is well-known in the literature of signal
processing.  For example in the mid-eighties, Santosa and Symes
proposed this rule to reconstruct spike trains from incomplete data
\cite{SantosaSymes}, see also \cite{DobsonSantosa}.  We would also
like to point out connections with total-variation approaches in the
literature of image processing \cite{TVDN,CRT} which are methods based
on the minimization of the $\ell_1$-norm of the discrete gradient.
Note that minimizing the $\ell_1$-norm is very different than standard
least squares (i.e. $\ell_2$) minimization procedures.  With
incomplete data, the least square approach would simply set to zero
the `unobserved' coefficients. Consider the Fourier case, for
instance. The least-squares solution would set to zero all the
unobserved frequencies so that the minimizer would have much smaller
energy than the original signal. As is well known, the minimizer would
also contain a lot of artifacts.

More recently, $\ell_1$-minimization perhaps best known under the name
of {\em Basis Pursuit}, has been proposed as a convex alternative to
the combinatorial norm $\ell_0$, which simply counts the number of
nonzero entries in a vector, for synthesizing signals as sparse
superpositions of waveforms \cite{BP}.  Interestingly, these methods
provided great practical success \cite{BP,ChenPhD} and were shown to
enjoy remarkable theoretical properties and to be closely related to
various kinds of uncertainty principles
\cite{DonohoHuo,DonohoElad,GribonvalNielsen,RandomBP}.

On the practical side, an $\ell_1$-norm minimization problem (for
real-valued signals) can be recast as a linear program (LP)
\cite{Bloomfield}. For example, \eqref{eq:P1} is equivalent to
minimizing $\sum_t u(t)$ subject to $F_\Omega g = F_\Omega f$ and
$-u(t) \le g(t) \le u(t)$ for all $t$. This is interesting since there
is a wide array of ever more effective computational strategies for
solving LPs.


\subsection{Applications}

In many applications of practical interest, we often wish to
reconstruct an object (a discrete signal, a discrete image and so on)
from incomplete samples and it is natural to ask how much one can hope
to recover. Actually, this work was motivated by the problem of
reconstructing biomedical images from vastly undersampled Fourier
data. Of special interest are problems in magnetic resonance (MR)
angiography but it is expected that our methodology and algorithms
will be suitable for other MR imagery, and to other acquisition
techniques, such as tomography. In MR angiography, however, we observe
few Fourier samples, and therefore if the images of interest are
compressible in some transform domain such as in the wavelet domain
for example, then $\ell_1$-based reconstructions might be especially
well-suited.

Another application of these ideas might be to view the
measurement/reconstruction procedure as a kind of lossy
encoder/decoder pair where the measurement process would play the role of
an encoder and the linear program $(P_1)$ that of a decoder. We
postpone this discussion to Section \ref{sec:universal}.

\subsection{Organization of the Paper}

This paper is roughly divided into three parts and is organized as
follows. The first part (Sections \ref{sec:l1} and \ref{sec:l2}),
shows how {\bf UUP} together with {\bf ERP} give our main result,
namely, Theorem \ref{general-lp-control}. In Section \ref{sec:l1}, we
establish that the solution to \eqref{eq:P1} is in some sense stable
in the $\ell_1$-norm, while Section \ref{sec:l2} introduces some
$\ell_2$-theory and proves our main result. In the second part
(Sections \ref{sec:eigen}, \ref{sec:low-level}, \ref{sec:ERP} and
\ref{sec:fourier}), we show that all three measurement ensembles obey
{\bf UUP} and {\bf ERP}. Section \ref{sec:eigen} studies singular
values of random matrices and shows that the {\bf UUP} holds for the
Gaussian and binary ensembles. Section \ref{sec:low-level} presents a
weaker {\bf ERP} which, in practice, is far easier to check. In
Section \ref{sec:ERP}, we prove that all three ensembles obey the {\bf
  ERP}. In the case of the Fourier ensemble, the strategy for proving
the {\bf UUP} is very different than for Gaussian and binary
measurements, and is presented in a separate Section
\ref{sec:fourier}.  Finally, we will argue in the third part of the
paper that one can think of the random measurement process as some
kind of universal encoder (Section \ref{sec:universal}) and briefly
discuss some of its very special properties. We conclude with a
discussion section (Section \ref{sec:discussion}) whose main purpose
is to outline further work and point out connections with the work of
others. The Appendix provides proofs of technical lemmas.

\section{Stability in the $\ell_1$-norm}
\label{sec:l1}

In this section, we establish $\ell_1$-properties of any minimizer to
the problem $(P_1)$, when the initial signal is mostly concentrated
(in an $\ell_1$ sense) on a small set.

\begin{lemma}
\label{teo:l1-stability}
Assume that the measurement matrix $F_\Omega$ obeys {\bf ERP}. We let
$f$ be a fixed signal of the form $f = f_0 + h$ where $f_0$ is a
signal supported on a set $T$ whose size obeys \eqref{eq:TK}. Then
with probability at least $1 - O(N^{-\rho/\alpha})$, any
$\ell_1$-minimizer \eqref{eq:P1} obeys
\begin{equation}
  \label{eq:l1-stability}
  \|f^\sharp \cdot 1_{T^c}\|_{\ell_1} \le 4 \,\|h \|_{\ell_1}. 
\end{equation}
\end{lemma}
\begin{proof} Observe that since $f$ is of course feasible for $(P_1)$, we
  immediately have
  \begin{equation}
    \label{eq:l1-trivial}
   \| f^\sharp \|_{\ell_1} \leq \|f\|_{\ell_1} \le \| f_0
  \|_{\ell_1} + \|h\|_{\ell_1}.
  \end{equation}
  Now because {\bf ERP} holds, one can construct---with the required
  probability---a function $P = F_\Omega^* V$ for some $V \in \ell_2(K)$ such that
  $P = \sgn(f_0)$ on $T$ and $|P(t)| \leq 1/2$ away from $T$. Observe the identity
$$
  \langle f^\sharp, P \rangle = \langle f^\sharp, F_\Omega^* V \rangle =
  \langle F_\Omega f^\sharp, V \rangle = \langle F_\Omega(f_0 + h), V
  \rangle = \langle f_0+h, F_\Omega^* V \rangle = \langle f_0 + h, P \rangle.
$$
Then on the one hand 
$$
\langle f^\sharp, P \rangle = \langle f_0, P \rangle + \<h, P\> \ge
\| f_0 \|_{\ell_1} - \|h\|_{\ell_1}, 
$$
while on the other hand, the bounds on $P$ give
  \begin{eqnarray*}
  |\langle f^\sharp, P \rangle| & \leq & \sum_T |f^\sharp(t) P(t)| +
  \sum_{T^c} |f^\sharp(t) P(t)|\\
  & \le & \sum_T |f^\sharp(t)| + \frac{1}{2} \sum_{T^c} |f^\sharp(t)|
  = \sum_{\ZZ_N} |f^\sharp(t)| - \frac{1}{2} \sum_{T^c}
  |f^\sharp(t)|.
  \end{eqnarray*}
  To conclude, we established that 
  $$
  \|f_0\|_{\ell_1} - \|h\|_{\ell_1} \le \|f^\sharp\|_{\ell_1} -
  \frac{1}{2} \|f^\sharp \, 1_{T^c}\|_{\ell_1},  
  $$
  and together with \eqref{eq:l1-trivial} proved that
  $$ 
 \|f^\sharp \, 1_{T^c}\|_{\ell_1} \le 4 \, \|h\|_{\ell_1},
 $$
 as claimed.
\end{proof}
This lemma says that any minimizer is approximately concentrated on
the same set as the signal $f$.  Indeed, suppose that $f$ obeys
\eqref{weak} and consider $T$ to be the set of largest values of
$|f|$.  Set $f_0 = f \cdot 1_T$.  Then the property \eqref{weak} gives 
$$
\|h\|_{\ell_1} = \|f \cdot 1_{T^c}\|_{\ell_1} \le C_p \cdot |T|^{1-1/p}
$$
for some constant $C_p$ only
depending on $p$, and therefore \eqref{eq:l1-stability} gives
\begin{equation}
  \label{eq:l1-bound}
  \|f^\sharp \cdot 1_{T^c}\|_{\ell_1} \le 4 C_p \cdot |T|^{1-1/p}.
\end{equation}
Thus, $f^\sharp$ puts `little mass' outside of the set $T$.

\begin{corollary}
\label{teo:linfty-stability}
Let $f^\sharp$ be any $\ell_1$-minimizer to the problem $(P_1)$ and
rearrange the entries of $f^\sharp$ in decreasing order of magnitude
$|f^\sharp|_{(1)} \ge |f^\sharp|_{(2)} \ge \ldots \ge
|f^\sharp|_{(N)}$.  Under the hypotheses of Lemma
\ref{teo:l1-stability}, the $m$th largest entry of $f^\sharp$ obeys
\begin{equation}
  \label{eq:linfty-stability}
  |f^\sharp|_{(m)} \le C_p \cdot |T|^{-1/p}, \quad \forall m > 2|T|. 
\end{equation}
\end{corollary}
\begin{proof}
  Suppose $T$ is the set of $|T|$ largest entries of $f$ as above so
  that $f^\sharp$ obeys \eqref{eq:l1-bound}. Denote by $E_m$ the set
  of the $m$-largest values of the function $f^\sharp$. Obviously
  $|E_m \cap T^c| \ge m - |T|$ and, therefore,
  $$
  \|f^\sharp\|_{\ell_1(E_m \cap T^c)} \ge (m- |T|) \cdot
  |f^\sharp|_{(m)} \ge |T| \cdot |f^\sharp|_{(m)}.$$
  The claim then follows from 
  $$
  \|f^\sharp\|_{\ell_1(E_m \cap T^c)} \le \|f^\sharp \, 1_{T^c}\|_{\ell_1} \le C
  \cdot |T|^{1-1/p}.
  $$
\end{proof}


\section{Stability in the $\ell_2$-norm}
\label{sec:l2}

\subsection{Extension lemma}

As essentially observed in \cite{CRT}, a matrix obeying
\eqref{eq:uup}--- think of it as a partial Fourier transform--- allows
to extend a function from a small set to all of $\ZZ_N$ while
constraining its Fourier transform to a fixed random set:
\begin{corollary}[Extension theorem]\label{choose}  
  Assume that $F_\Omega$ is a matrix obeying the uniform uncertainty
  principle {\bf UUP}.  Then with probability at least $1 -
  O(N^{-\rho/\alpha})$ the following statement holds: for all sets $T
  \subset \ZZ_N$ obeying the bound \eqref{eq:uup} and all functions $f
  \in \ell_2(T)$, there exists $f^{\operatorname{ext}} \in \ell_2(\ZZ_N)$ which
  \begin{itemize}
  \item agrees  with $f$ on $T$ ($R_T f^{\operatorname{ext}} = f$), 
  \item belongs to the column space of $F_\Omega^*$ (i.e. $f^{\operatorname{ext}} =
    F_\Omega^* V$ for some $V \in \ell_2(\Omega)$),
  \item and furthermore, we have the $\ell_2$ estimates
\begin{equation}\label{F-disjoint}
 \| f^{\operatorname{ext}} \|_{\ell_2(E)} \leq 
C \left(1 + \frac{|E|}{\alpha 
K/\lambda}\right)^{1/2} \|f\|_{\ell_2(T)}
\end{equation}
valid for all $E \subseteq \ZZ_N$.  
  \end{itemize}
\end{corollary}
\begin{proof}  We may assume that we are on an event such 
  that the conclusions of {\bf UUP} hold.  In particular, from
  \eqref{eq:uup}, the operator $(F_{\Omega T}^* F_{\Omega T})$ is
  invertible and the inverse obeys 
\begin{equation}\label{neumann-consequence}
 \|(F_{\Omega T}^* F_{\Omega T})^{-1}\| \leq 2 N/K,
 \end{equation}
 where $\|\cdot\|$ is the operator norm induced by the $\ell_2$ norm.
 In the remainder of this paper and unless specified otherwise $\|A\|$
 will always be the operator norm of $\|A\| := \sup_{\|x\|_{\ell_2} =
   1} \|Ax\|_{\ell_2}$.  We now set $f^{\operatorname{ext}}$ as 
$$
f^{\operatorname{ext}} := F_\Omega^* F_{\Omega T} \, (F_{\Omega T}^* F_{\Omega T})^{-1} \, f.
$$
By construction, $f^{\operatorname{ext}}$ agrees with $f$ on $T$, and is in the
column space of $F_\Omega^*$.  Now we prove \eqref{F-disjoint}.  It suffices
to do so when $|E| \leq \alpha K/\lambda$, since the general claim
then follows by decomposing larger $E$'s into smaller sets and then
square-summing.  But from \eqref{eq:uup}, we see that $F_{\Omega T}$
and $F^*_{\Omega E}$ have operator norms of size at most
$\sqrt{3K/2N}$, and the claim follows by composing these facts with
\eqref{neumann-consequence}.
\end{proof}

\subsection{Proof of Theorem \ref{general-lp-control}}\label{gen-sec}
\newcommand{\qed}{\quad\hbox{\vrule width 4pt height 6pt depth 1.5pt}}

Let $T_0$ (resp.~$T_1$) be the set of the $S$-largest values of $|f|$
(resp.~$|f^\sharp|$) and put $T = T_0 \cup T_1$. By construction, $S
\le |T| \le 2S$ and we assume that $|T|$ obeys the condition
\eqref{eq:TK}. Now observe that by construction of $T$, a consequence of 
Lemma \ref{teo:l1-stability} is that
\begin{equation}\label{fsharp-l1}
 \| f - f^\sharp \|_{\ell_1(T^c)} \leq \|f\|_{\ell_1(T_0^c)} + 
\|f^\sharp\|_{\ell_1(T_0^c)} \le C_p \cdot |T|^{1-1/p}.
 \end{equation}
Furthermore, it follows from our assumption about $f$ and 
\eqref{eq:linfty-stability} that 
\begin{equation}\label{fsharp-linfty}
  \| f - f^\sharp\|_{\ell_\infty(T^c)} \leq \| f \|_{\ell_\infty(T_0^c)} + \| f^\sharp\|_{\ell_\infty(T_1^c)} \leq C \cdot |T|^{-1/p}.
 \end{equation}
By interpolation, these last two inequalities give 
\begin{equation}\label{fsharp-l2}
 \| f - f^\sharp\|_{\ell_2(T^c)} \leq C \cdot |T|^{1/2-1/p},
 \end{equation}
and it remains to prove that the same bound holds over the set $T$.

In order to prove this fact,
Corollary \ref{choose} assures us that one can find a
function of the form $g = F_\Omega^* V$ which matches $h$ on $T$ and
with the following property:
\begin{equation}\label{g-josh}
\sum_{t \in E} |g(t)|^2 \leq C \sum_{t \in T} |f(t) - f^\sharp(t)|^2, 
\end{equation}
for all sets $E$ of cardinality $O(K/\lambda)$ that are disjoint from
$T$. Here and in the rest of the proof, the constants $C$ are allowed
to depend on $\alpha$. From the representation $g = F_\Omega^* V$ and
the constraint $F_\Omega f = F_\Omega f^\sharp$ (from \eqref{eq:P1}),
we have
\[
\langle f - f^\sharp, g \rangle = \langle f - f^\sharp, F_\Omega^* V 
\rangle = \langle F_\Omega f - F_\Omega f^\sharp, V \rangle = 0,
\]
and hence
$$ \sum_{t \in\ZZ_N} (f - f^\sharp)(t) \overline{g(t)} = 0.
$$
Splitting into $T$ and $T^c$, we obtain
\begin{equation}\label{F-amplitude}
 \sum_{t \in T} |f - f^\sharp|^2(t) = - \sum_{t \in T^c}
(f - f^\sharp)(t) \overline{g(t)}.
\end{equation}

We will use \eqref{F-amplitude} to show that the left-hand side must
be small since \eqref{fsharp-l2} and \eqref{g-josh} assert that the
right-hand side is not very large.  Enumerate $T^c$ as $n_1, n_2,
\ldots, n_{N-|T|}$ in decreasing order of magnitude of
$|f-f^\sharp|$. We then group these into adjacent blocks $B_J$ of size
$|T|$ (except perhaps for the last one) $B_J := \{n_j, J |T| < j \le
(J + 1) |T|\}$, $J = 0, 1, \ldots$.  From \eqref{g-josh} and
Cauchy-Schwarz, we have
\begin{equation}
  \label{eq:intermediate}
\left| \sum_{j \in B_J}
(f - f^\sharp)(n_j) \overline{g(n_j)} \right| 
\le C \cdot \| f - f^\sharp \|_{\ell_2(T)} \cdot I_J,
\end{equation}
where 
\[
I_J := \sqrt{\sum_{j = J |T| + 1}^{(J+1) |T|} |(f -
  f^\sharp)(n_j)|^2}.
\]
Because we are enumerating the values of $n_j$ in decreasing order, we
have $I_0 \le |T|^{1/2} \cdot |f - f^\sharp|(n_1) \le C \cdot
|T|^{1/2-1/p}$ while for $J \ge 1$ we have
\[
I_J \le |T|^{1/2} \cdot |f - f^\sharp|(n_{J |T| + 1}) \le |T|^{1/2} \cdot
|T|^{-1}\cdot \|f - f^\sharp\|_{\ell_1(B_{J-1})},
\]
In other words, 
\[
\sum_J I_J \le I_0 + \sum_{J \ge 1} I_J \le C\cdot |T|^{-r} +
|T|^{-1/2} \cdot \|f - f^\sharp\|_{\ell_1(T^c)} 
\]
and, therefore, it follows from \eqref{fsharp-l1} that the summation of 
 the inequality \eqref{eq:intermediate} over the blocks $B_J$ gives
$$
\left| \sum_{t \in T^c} (f - f^\sharp)(t) \overline{g(t)} \right| \le
C \cdot |T|^{-r} \cdot \|f - f^\sharp \|_{\ell_2(T)}.
$$
Inserting this back into \eqref{F-amplitude}, we established 
$$ 
\|f - f^\sharp \|_{\ell_2(T)} \le C \cdot |T|^{-r}.
$$
This concludes the proof of Theorem \ref{general-lp-control}. $\qed$

Note that by Cauchy-Schwarz, it follows from the proof of our Theorem
that
$$
\|f - f^\sharp \|_{\ell_1(T)} \le C \cdot |T|^{1-1/p}, 
$$
and, therefore, owing to \eqref{fsharp-l1}, we also proved an
$\ell_1$ stability estimate
\begin{equation}
  \label{eq:l1-stability2}
  \|f - f^\sharp \|_{\ell_1} \le C \cdot |T|^{1-1/p}.
\end{equation}

Had we assumed that $f$ belonged to the weak-$\ell_1$ ball when $p =
1$, the right-hand side of \eqref{fsharp-l1} would read $C_1 \log
(N/|T|)$ instead of just $C_1$. This is the reason why we required
$\ell_1$ in the hypothesis of Theorem \ref{general-lp-control} and
showed that we also have a near-optimal signal recovery result for the
unit ball of $\ell_1$ with no additional losses (logarithmic or
otherwise).

\subsection{Uniqueness of the minimizer for the Gaussian ensemble}

The claim that the minimizer $f^\sharp$ is unique with probability 1,
for Gaussian measurements, can be easily established as follows.  The
claim is trivial for $f \equiv 0$ so we may assume $f$ is not
identically zero.  Then $F_\Omega f$ is almost surely non-zero.
Furthermore, if one considers each of the (finitely many) facets of
the unit ball of $\ell_1(\ZZ_N)$, we see that with probability 1 the
random Gaussian matrix $F_\Omega$ has maximal rank on each of these
facets (i.e. the image of each facet under $F_\Omega$ has dimension
equal to either $K$ or the dimension of the facet, whichever is
smaller).  From this we see that every point on the boundary of the
image of the unit $\ell_1$-ball under $F_\Omega$ arises from a unique
point on that ball.  Similarly for non-zero dilates of this ball.
Thus the solution to the problem \eqref{eq:P1} is unique as claimed.

We remark that the question of establishing uniqueness with high
probability for discretely randomized ensembles such as the binary and
Fourier ensembles discussed below is an interesting one, but one which
we will not pursue here.

\section{Eigenvalues of random matrices}
\label{sec:eigen}

In this section, we show that all three ensembles obey the uniform
uncertainty principle {\bf UUP}. 

\subsection{The Gaussian ensemble}

Let $X$ be an $n$ by $p$ matrix with $p \le n$ and with i.i.d.~entries
sampled from the normal distribution with mean zero and variance
$1/n$.  We are interested in the singular values of $X$ or the
eigenvalues of $X^* X$. A famous result due to Marchenko and Pastur
\cite{MarchenkoPastur} states that the eigenvalues of $X^* X$ have a
deterministic limit distribution supported by the interval $[(1 -
\sqrt{c})^2, (1+ \sqrt{c})^2]$ as $n,p \goto \infty$, with $p/n \goto
c < 1$. In fact, results from \cite{Silverstein} show that the
smallest (resp.~largest) eigenvalue converges a.s.~to $(1-
\sqrt{c})^2$ (resp.~$(1 + \sqrt{c})^2$). In other words, the smallest
singular value of $X/\sqrt{n}$ converges a.s.~to $1 - \sqrt{c}$ and
the largest to $1 + \sqrt{c}$.  In addition, there are remarkably fine
statements concerning the speed of the convergence of the largest
singular value \cite{IainWishart}.

To derive the {\bf UUP}, we need a result about the concentration of
the extreme singular values of a Gaussian matrix, and we borrow a most
elegant estimate due to Davidson and Szarek \cite{Szarek2}. We let
$\lambda_1(X) \le \ldots \le \lambda_p(X)$ be the ordered list of the
singular values of $X$. Then in \cite{Szarek2}, the authors prove that
\begin{align}
\P\left(\lambda_p(X) > 1 + \sqrt{p/n} + r\right)  & \le e^{- nr^2/2}, 
\label{eq:concentration4a}\\
\P\left(\lambda_1(X) < 1 - \sqrt{p/n} - r\right)  & \le e^{- nr^2/2}.
\label{eq:concentration4b}
\end{align}
Such inequalities about the concentration of the largest and smallest
singular values of Gaussian matrices have been known for at least a
decade or so. Estimates similar to
\eqref{eq:concentration4a}-\eqref{eq:concentration4b} may be found in
the work of Szarek \cite{Szarek1}, see also Ledoux \cite{LedouxAMS}.

\begin{lemma}
\label{teo:GaussUUP}
  The Gaussian ensemble obeys the uniform uncertainty principle ({\bf UUP}) with
  oversampling factor $\lambda = \log N$. 
\end{lemma}
\begin{proof} Fix $K \geq \log N$ and let $\Omega := \{1,\ldots,K\}$.
  Let $T$ be a fixed index set and define the event $E_T$ as
\[
E_T := \{\lambda_{\min} (F_{\Omega T}^* F_{\Omega T}) < K/2N\} \,\,\, \cup \,\,\,
 \{\lambda_{\max}(F_{\Omega T}^* F_{\Omega T}) >  3K/2N\}.
\]
Since the entries of $F_{\Omega T}$ are i.i.d. $N(0,1/N)$, it follows
  from \eqref{eq:concentration4a}-\eqref{eq:concentration4b} by a
  simple renormalization that for each $|T| \le K/16$,
\[
\P(E_T) \le 2 e^{-cK}
\]
where one can choose $c = 1/32$ by selecting $r = 1/4$ in
\eqref{eq:concentration4a}-\eqref{eq:concentration4b}. We now examine
the tightness of the spectrum over all sets $T \in {\cal T}_m := \{T :
|T| \le m\}$ where we assume that $m$ is less than $N/2$.  We have
\[
\P\left(\cup_{{\cal T}_m} E_T \right) \le 2 e^{-c K} \cdot |{\cal
  T}_m| = 2 e^{-c K} \cdot \sum_{k = 1}^{m} \binom{N}{k} \le 2 e^{-c
  K} \cdot m \, \binom{N}{m}.
\]
We now use the well-known bound on the binomial coefficient 
\[
\binom{N}{m} \le e^{N H(m/N)},
\]
where for $0 < q < 1$, $H$ is the binary entropy function
\[
H(q) := - q\log q - (1-q) \log(1-q). 
\]
The inequality $-(1-q)
\log(1-q) \le q$ shows that $-(1-m/N)\log(1-m/N) \le m/N$ and thus 
\[
m \cdot \binom{N}{m} = e^{m\log(N/m) + m + \log m}.
\]
Whence, 
\[ 
\log \P\left(\cup_{{\cal T}_m} E_T \right) \le \log 2 - c K +
m(\log(N/m) + 1 + m^{-1} \log m) \le \log 2 - \rho K
\]
provided that $m (\log(N/m) + 1 + m^{-1} \log m) \le (c - \rho) K$,
which is what we needed to establish. (We need to assume that $K \ge
(c-\rho)^{-1}(1+\log N)$ for the claim not to be vacuous.)  Note that
we proved more than what we claimed since the {\bf UUP} holds for an
oversampling factor proportional to $\log N/K$.
\end{proof}

\subsection{The binary ensemble}

The analysis is more complicated in the case where the matrix $X$ is an
$n$ by $p$ array with i.i.d. symmetric Bernoulli entries taking on
values in $\{-1/\sqrt{n}, 1/\sqrt{n}\}$. To study the concentration of
the largest singular values of $X$, we follow an approach proposed by
Ledoux \cite{LedouxAMS} which makes a simple use of the concentration
property, see also \cite{ElKarouiPhD}. 

As before, we let $\lambda_p(X)$ be the mapping that associates to a
matrix $X$ its largest singular values. Equip $\RR^{np}$ with the
Frobenius norm
\[
\|X\|^2_{F} := \tr(X^* X) = \sum_{i,j = 1}^p |X_{ij}|^2
\]
(the Euclidean norm over $\RR^{np}$). Then the mapping $\lambda_p$ is
convex and 1-Lipschitz in the sense that
\[
|\lambda_p(X) - \lambda_p(X')| \le \|X - X'\|_{F}
\]
for all pairs $(X, X')$ of $n$ by $p$ matrices. A classical
application of the concentration inequality for binary measures
\cite{LedouxAMS} 
then gives
\begin{equation}
\label{eq:concentration}
\P\left(\lambda_{p}(X) - m(\lambda_{p}(X)) \ge r\right) \le e^{- nr^2/16},   
\end{equation}
 $m(\lambda_{p}(X))$ is either the mean or the median of
$\lambda_{p}(X)$. Now the singular values still exhibit the same
behavior; that is $\lambda_{\min}(X/\sqrt{n})$ and
$\lambda_{\max}(X/\sqrt{n})$ converge a.s. to $1 - \sqrt{c}$ and $1 +
\sqrt{c}$ respectively, as $n, p \goto \infty$ with $p/n \goto c$
\cite{BaiYin}.  As a consequence, for each $\epsilon_0$ and $n$
sufficiently large, one can show that the medians belong to the fixed
interval $[1-\sqrt{p/n} - \epsilon_0, 1 + \sqrt{p/n} + \epsilon_0]$
which gives
\begin{equation}
  \label{eq:concentration4c}
\P\left(\lambda_{p}(X) >  1 + \sqrt{p/n} + \epsilon_0 + r\right) \le 
e^{- nr^2/16}. 
\end{equation}
This is a fairly well-established result \cite{ElKarouiPhD}.

The problem is that this method does not apply to the minimum singular
value which is 1-Lipshitz but not convex. Fortunately, Litvak, Pajor,
Rudelson and Tomczak-Jaegermann \cite{LitvakPajor}[Theorem 3.1] have
recently announced a result which gives exponential concentration for
the lowest singular value. They proved that whenever $n \ge
(1+\delta)p$ where $\delta$ is greater than a small constant, 
\begin{equation}
  \label{eq:concentration4d}
  \P\left(\lambda_{1}(X) >  c_1\right) \le 
    e^{- c_2 n}, 
\end{equation}
where $c_1$ and $c_2$ are universal positive constants.

Just as \eqref{eq:concentration4a}-\eqref{eq:concentration4b} implied
the uniform uncertainty principle {\bf UUP} for Gaussian matrices,
\eqref{eq:concentration4c}-\eqref{eq:concentration4c} gives the same
conclusion for the binary ensemble with the proviso that the condition
about the lowest singular value reads $\lambda_{\min}(F_{\Omega T}^*
F_{\Omega T}) > c_1 K/N$; i.e., $c_1$ substitutes 1/2 (recall the
remark following the definition of the {\bf UUP}).
\begin{lemma}
\label{teo:binaryUUP}
The binary ensemble obeys the uniform uncertainty principle ({\bf UUP})
with oversampling factor $\lambda = \log N$.
\end{lemma}
The proof is of course identical to that of Lemma
\ref{teo:GaussUUP}. If we define $E_T$ as
\[
E_T := \{\lambda_{\min}(F_{\Omega T}^* F_{\Omega T}) < c_1 K/N\} \,\,\, \cup \,\,\,
 \{\lambda_{\max}(F_{\Omega T}^* F_{\Omega T}) >  3K/2N\}, 
\]
we have $\P(E_T) \le 2 e^{-cK}$ for some constant $c > 0$. The rest of
the proof is as before.

\subsection{The Fourier ensemble}

The analysis for the Fourier ensemble is much more delicate than that
for the Gaussian and binary cases, in particular requiring entropy
arguments as used for instance by Bourgain \cite{bourgain-lambda},
\cite{bourgain-halasz}.  We prove the following lemma in the separate
Section \ref{sec:fourier}.
\begin{lemma}
\label{teo:FourierUUP}
The Fourier ensemble obeys the uniform uncertainty principle {\bf UUP} with
oversampling factor $\lambda = (\log N)^6$.
\end{lemma}

The exponent of 6 can almost certainly be lowered\footnote{Note added
  in proof: since the submission of this paper, Rudelson and
  Vershynin, in a very recent piece of work \cite{rv}, have improved
  the oversampling factor to $(\log N)^4$.}, but we will not attempt
to seek the optimal exponent of $\log N$ here.

\section{Generic signals and the weak ERP}
\label{sec:low-level}

In some cases, it might be difficult to prove that the exact
reconstruction principle {\bf ERP} holds, and it is interesting to
observe that {\bf UUP} actually implies {\bf ERP} for `generic' sign
functions $\sigma = \pm 1$ supported on a small set $T$.  More
precisely, if we fix $T$ and define $\sigma$ to be supported on $T$
with the i.i.d. Bernoulli distribution (independently of $F_\Omega$),
thus
\[
\P(\sigma(t)  = \pm 1) = 1/2, \quad \text{for all } t \in T.
\]
then we shall construct a $P$ obeying the conditions (i)-(iii) in the
definition of {\bf ERP}.  Indeed, we shall construct $P$ explicitly
as
\begin{equation}
  \label{eq:l2-P}
  P := F_\Omega^* F_{\Omega T} (F_{\Omega T}^* F_{\Omega T})^{-1} R_T \sigma; 
\end{equation}
one can view this choice of $P = F_\Omega^* V$ as the unique solution
to (i) and (ii) which minimizes the $\ell_2$ norm of $V$, and can thus
be viewed as a kind of least-squares extension of $\sigma$ using the
rows of $F_\Omega$.

It is immediate to check that $P$ obeys (i) and (ii) above.  Indeed,
the restriction of $P$ to $T$ is given by
\[
R_T P = R_T F_\Omega^* F_{\Omega T} (F_{\Omega T}^* F_{\Omega T})^{-1}
R_T \sigma = F_{\Omega T}^* F_{\Omega T} (F_{\Omega T}^* F_{\Omega
  T})^{-1} R_T \sigma = R_T \sigma
\] 
and, therefore, (i) is verified. Further, it follows from the
definition that $P$ is a linear combination of the columns of
$F_\Omega^*$ and thus, (ii) holds. Therefore, we only need to check
that for all $t \in T^c$, $|P(t)| \le \frac{1}{2}$ with sufficiently
high probability.  In order to do this, we rewrite $P(t)$ as
\[
P(t) = \<W_t, R_T \sigma\>, 
\]
where for each $t \in T^c$, $W_t$ is the $|T|$ dimensional vector
\[
W_t := (F_{\Omega T}^* F_{\Omega T})^{-1} \, F_{\Omega T}^* F_t 
\]
and $F_t$ is the $t$-{th} column of $F_\Omega$. We now introduce
another condition which is far easier to check than {\bf ERP}.

{\bf WERP (Weak ERP).} {\em We say the the measurement process obeys
  the weak {\bf ERP}, if for each fixed $T$ obeying \eqref{eq:TK} and
  any $0 < \gamma \leq 1$,
  $F_\Omega$ obeys
\begin{equation}
  \label{eq:easy}
  \| F_{\Omega T}^* F_t\|_{\ell_2} \le \gamma \cdot \sqrt{K \, |T|}/N \quad 
\text{for all } t \in T^c 
\end{equation}
with probability at least $1 - O(N^{-\rho/\gamma})$ for some fixed
positive constant $\rho > 0$.}

For example, it is an easy exercise in large deviation theory to show
that {\bf WERP} holds for Gaussian and binary measurements. One can
also check that {\bf WERP} holds for random frequency samples. We
omit the proof of these facts, however, since we will show the
stronger version, namely, {\bf ERP} in all three cases.  Instead, we
would like to emphasize that {\bf UUP} together with {\bf WERP}
actually imply {\bf ERP} for most sign patterns $\sigma$.

We begin by recalling the classical Hoeffding inequality: let
$X_1,\ldots,X_{N} = \pm 1$ be independent symmetric Bernoulli random variables
and consider the sum $S = \sum_{j = 1}^N a_j X_j$. Then
\begin{equation}
  \label{eq:Hoeffding}
  P(|S| \le \lambda) \le 2 
\exp\left(-\frac{\lambda^2}{2\|a\|_{\ell_2}^2}\right). 
\end{equation}
Suppose now that the $\sigma(t)$'s are independent Bernoulli, and
independent from $F$. Then \eqref{eq:Hoeffding} gives 
\[
\P\left(|P(t)| > \lambda \,\, | \,\, \|W_t\|_{\ell_2} = \rho\right) \le  2 
\exp\left(-\frac{\lambda^2}{2\rho^2}\right). 
\]
If we now assume that both {\bf UUP} and {\bf WERP} hold, then for any $0 < \gamma \leq 1$ we have
\[
\|W_t\|_{\ell_2} \le \|(F_{\Omega T}^* F_{\Omega T})^{-1}\| \cdot
\| F_{\Omega T}^* F_t\|_{\ell_2} \le 2\gamma \cdot
\sqrt{\frac{|T|}{K}}.
\]
with probability at least $1 - O(N^{-\rho/\gamma})$. This shows that 
\begin{eqnarray*}
  \P(|P(t)| > \lambda) & \le & 
\P\left(|P(t)| > \lambda \,\, \vert \,\,
  \|W_t\|_{\ell_2} \le 2\gamma \cdot \sqrt{|T|/K}\right) + 
\P\left(\|W_t\|_{\ell_2} > 2\gamma \cdot \sqrt{|T|/K}\right)\\
& \le &  2 \exp\left(-\frac{\lambda^2}{8\gamma^2}\, 
\frac{K}{|T|} \right) +  O(N^{-\rho/\gamma}) = 2 
e^{-\rho_0 K/|T|} +  O(N^{-\rho/\gamma}). 
\end{eqnarray*}
Hence, if $|T| \le \alpha K/\log N$, then
\[
\P\left( \sup_{t \in T^c} |P(t)| > 1/2 \right) 
\le O(N \cdot N^{-\rho_0/\alpha}) + O(N^{-\rho/\alpha}).  
\]
Therefore, if $\alpha$ is chosen small enough, then for some small
$\rho' > 0$
\[
\P\left( \sup_{t \in T^c} |P(t)| > 1/2 \right) = 
O(N^{-\rho'/\alpha}).
\]
In other words, {\bf ERP} holds for most sign patterns $\sigma$. That
is, if one is only interested in providing good reconstruction to
nearly all signals (but not all) in the sense discussed above, then it
is actually sufficient to check that both conditions {\bf UUP} and
{\bf WERP} are valid.

\section{About the exact reconstruction principle}
\label{sec:ERP}

In this section, we show that all the three ensembles obey the exact
reconstruction principle {\bf ERP}.

\subsection{The Gaussian ensemble}

To show that there is function $P$ obeying the conditions (i)-(iii) in
the definition of {\bf ERP}, we take an approach that resembles that
of Section \ref{sec:low-level}, and establish that $P$ defined as in
\eqref{eq:l2-P}, 
\[
P := F_{\Omega}^* F_{\Omega T} (F_{\Omega T}^* F_{\Omega T})^{-1} R_T \sigma
\]
obeys the three conditions (i)-(iii).

We already argued that $P$ obeys (i) and (ii). Put $P^c = R_{T^c} P$ 
to be the restriction of $P$ to $T^c$. We need to show that 
\[
\sup_{T^c} |P^c(t)| \le 1/2
\]
with high probability. Begin by factorizing $P^c$ as 
\[
P^c =  F^*_{\Omega T^c} V, \hbox{ where } V := (F_{\Omega T}^* F_{\Omega T})^{-1} R_T \sigma. 
\]
The crucial observation is that the random matrix $F^*_{\Omega T^c}$ and 
the random variable $V$ are independent
since they are functions of disjoint sets of independent variables.  
\begin{proposition}
\label{teo:simple}
Conditional on $V$, the components of $P^c$ are i.i.d. Gaussian with 
\[
P^c(t) \sim N(0, \|V\|^2/N). 
\]
\end{proposition}
\begin{proof}
Suppose $V$ is fixed. By definition,
\[
P^c(t) = \frac{1}{\sqrt{N}} \sum_{k \in \Omega} X_{k,t} V_k  
\]
and, therefore, it follows from the independence between the
$X_{k,t}$'s and $V$ for each $t \in T^c$ that the conditional
distribution of $P^c(t)$ is normal with mean 0 and variance $\|V\|_{\ell_2}^2/N$.
The independence between the components of $P^c$ is a simple
consequence of the independence between the columns of $F$.
\end{proof}

\begin{lemma}  Let $\alpha > 0$ be sufficiently small, and suppose that $|T|$ is chosen as in \eqref{eq:TK} so that {\bf UUP}
  holds.  The components of $P^c(t)$ obey
\[
\P(|P^c(t)| > \lambda) \le \P\left(|Z| > \lambda \cdot
  \sqrt{K/6|T|}\right) + O(N^{-\rho/\alpha}). 
\]
where $Z \sim N(0,1)$ is a standard normal random variable.
\end{lemma}
\begin{proof}
Observe that 
\begin{equation*}
\|V\|_{\ell_2}  \le  \|F_T\| \cdot \|(F_T^* F_T)^{-1}\| \cdot \|R_T \sigma\|_{\ell_2}.
\end{equation*}
On the event such that the conclusions of {\bf UUP} holds, $\|F_T\|
\le \sqrt{3K/2N}$ and also $\|(F_T^* F_T)^{-1}\| \le 2N/K$.  Since
$\|R_T \sigma\|_{\ell_2} = \sqrt{|T|}$, this gives
\[
\|V\|_{\ell_2} \le \sqrt{\frac{6N \cdot |T|}{K}}. 
\]
Therefore, 
\[
\P(|P^c(t)| > \lambda) \le \P\left(|P^c(t)| > \lambda \,\, \vert \,\,
  \|V\|_{\ell_2} \le \sqrt{6N |T|/K}\right) + \P\left(\|V\|_{\ell_2} > \sqrt{6N
    |T|/K}\right).
\]
The first term is bounded by Proposition \ref{teo:simple} while the
second is bounded via the uniform uncertainty principle {\bf UUP}. This
establishes the lemma.
\end{proof}

The previous lemma showed that for $|T| \le \alpha \cdot K/\log N$,
\[
\P\left( \sup_{t \in T^c} |P^c(t)| > 1/2 \right) 
\le N \cdot \P\left(|Z| > \frac{1}{2} \, \sqrt{\frac{\log N}{6 \alpha}}\right) + 
O(N^{-\rho/\alpha}) \le 2N^{1 - 1/(48 \alpha)} + O(N^{-\rho/\alpha}).
\]
Therefore, if $\alpha$ is chosen small enough, then for some small
$\rho' > 0$
\[
\P\left( \sup_{t \in T^c} |P^c(t)| > 1/2 \right) = 
O(N^{-\rho'/\alpha}).
\]
In short, we proved:
\begin{lemma}
\label{teo:GaussERP}
The Gaussian ensemble obeys the exact reconstruction principle {\bf ERP} with
oversampling factor $\lambda = \log N$.
\end{lemma}

\subsection{The binary ensemble}

The strategy in the case where the entries of $F$ are independent
Bernoulli variables is nearly identical and we only discuss the main
differences. Define $P$ and $V$ as above; obviously, $F^*_{\Omega T^c}$ and
$V$ are still independent.
\begin{proposition}
\label{teo:simple-binary}
Conditional on $V$, the components of $P^c$ are independent and obey 
\[
\P(|P^c(t)| > \lambda \,\, | \,\, V) \le 2 \exp\left(-\frac{\lambda^2
    N}{2\|V\|_{\ell_2}^2}\right).
\]
\end{proposition}
\begin{proof}
  The conditional independence of the components is as before. As far
  as the tail-bound is concerned, we observe that $P^c(t)$ is a
  weighted sum of independent Bernoulli variables and the claim
  follows from the Hoeffding inequality \eqref{eq:Hoeffding}.
\end{proof}

The rest of the argument is as before. If $|T|$ is selected as in
\eqref{eq:TK} such that  {\bf UUP} holds, one has
\[
\P(|P^c(t)| > 1/2) \le 2 N^{-1/48 \alpha} + O(N^{-\rho/\alpha}).
\]
And, of course, identical calculations now give
\begin{lemma}
\label{teo:binaryERP}
The binary ensemble obeys the exact reconstruction principle {\bf ERP} with
oversampling factor $\lambda = \log N$.
\end{lemma}

\subsection{The Fourier ensemble}

It turns out that the exact reconstruction principle also holds for
the Fourier ensemble although the argument is considerably more
involved \cite{CRT}. We do not reproduce the proof here but merely
indicate the strategy for proving that $P$ (defined as before) also
obeys the desired bound on the complement of $T$ with sufficiently
high probability.  We first remark that $|\Omega|$ is concentrated
around $K$.  To see this, recall the Bernstein's inequality
\cite{MassartSharp} which states that if $X_1, \ldots, X_m$ are
independent random variables with mean-zero and obeying $|X_i| \le c$,
then
\begin{equation}
  \label{eq:bernstein}
\P\left(\left| \sum_{i = 1}^N X_i \right| > \lambda \right) \le
2\exp\left(-\frac{\lambda^2}{2\sigma^2 + 2c\lambda/3}\right),  
\end{equation}
where $\sigma^2 = \sum_{i = 1}^m \mbox{\rm Var}(X_i)$.  Specializing
this inequality gives the following lemma which we shall need later
in this paper.
\begin{lemma}
\label{ms-lemma}
Fix $\tau \in (0,1)$ and let $I_k \in \{0,1\}$ be an i.i.d. sequence
of random variables obeying $\P(I_k = 1) = \tau$.  Let $a \in
\ell_2(\Z_N)$ be arbitrary, and set $\sigma^2 := \tau (1-\tau) \|a
\|_{\ell_2}^2$.  Then
  $$
  \P( |\sum_{k \in \Z_N} (I_k - \tau) a(k)| > \lambda) \leq 4
  \exp\left(-\frac{\lambda^2}{4 (\sigma^2 +
      \lambda \|a\|_{\infty}/3\sqrt{2})}\right).$$
\end{lemma}
\begin{proof}  Letting $S$ be the sum $\sum_{k \in \Z_N} (I_k - \tau) a(k)$, 
  the proof follows from \eqref{eq:bernstein} by simply stating that
  $P(|S| > \lambda)$ is bounded above by the sum $P(|S_1| >
  \lambda/\sqrt{2}) + P(|S_2| > \lambda/\sqrt{2})$, where $S_1$ and
  $S_2$ are the real and imaginary parts of $S$ respectively.
\end{proof}
Thus the bound on the quantity $|\sum_{k \in \Z_N} (I_k - \tau) a(k)|$
exhibits a Gaussian-type behavior at small thresholds $\lambda$, and an
exponential-type behavior at large thresholds.

Recall that $K = \E(|\Omega|) = N \tau$.  Applying Lemma
\ref{ms-lemma} with $a \equiv 1$ (so $\sigma^2 = N \tau(1-\tau)$), we
have that $P(K/2 \le |\Omega| \le 2K)$ with probability
$O(N^{-\rho/\alpha})$ provided that $K \ge C \, \alpha^{-1} \, \log
N$, which we will assume as the claim is vacuous otherwise.  In the
sequel, we assume that we are on an event $\{K/2 \le |\Omega| \le
2K\}$.

Decompose $F_{\Omega T}^* F_{\Omega T}$ as
\[
F_{\Omega T}^* F_{\Omega T} = \frac{|\Omega|}{N}(I - H), 
\]
where $H$ is the matrix defined by $H(t,t') = |\Omega|^{-1} \sum_{k
  \in \Omega} e^{i 2\pi k (t - t')}$ if $t \neq t'$ and $0$ otherwise.
We then expand the inverse as a truncated Neumann series
\[
(F_{\Omega T}^* F_{\Omega T})^{-1} = \frac{N}{|\Omega|}(I + H + \ldots
+ H^{n} + E),
\]
where $E$ is small remainder term. This allows to express $P^c$ as
\[
P^c = \frac{N}{|\Omega|} \cdot F_{\Omega T^c}^* F_{\Omega T} \cdot (I
+ H + \ldots + H^{n} + E),
\]
and one can derive bounds on each individual terms so that the sum
obeys the desired property. By pursuing this strategy, the following
claim was proved in \cite{CRT}.
\begin{lemma}
\label{teo:FourierERP}
The Fourier ensemble obeys the exact reconstruction principle {\bf ERP} with
oversampling factor $\lambda = \log N$.
\end{lemma}

\section{Uniform Uncertainty Principles for the Fourier Ensemble}
\label{sec:fourier}
\label{bourgain-sec}

In this section, we prove Lemma \ref{teo:FourierUUP}.  The ideas here
are inspired by an entropy argument sketched in
\cite{bourgain-halasz}, as well as by related arguments in
\cite{bourgain-lambda}, \cite{ms}.  These methods have since become
standard in the high-dimensional geometry literature, but we shall
give a mostly self-contained presentation here.

We remark that the arguments in this section (and those in the
Appendix) do not use any algebraic properties of the Fourier transform
other than the Plancherel identity and the fact that the maximum entry
of the Fourier matrix is bounded by $1/\sqrt{N}$. Indeed a simple
modification of the arguments we give below also gives the UUP for
randomly sampled rows of orthonormal matrices, see also
\cite{bourgain-lambda} and \cite{CR} for further discussion of this
issue. Suppose that $\sup_{i,j} |U_{ij}| \le \mu$ and let $U_\Omega$
be the matrix obtained by randomly selected rows. Then the UUP holds
for
\[
|T| \le C \cdot \frac{\E |\Omega|}{\mu^2 \, \log^6 N}.
\] 
In the case where one observes a few coefficients in the basis $\Phi$
when the signal is sparse in another basis $\Psi$, $\mu = \sqrt{N} \,
\sup_{i,j} |\<\phi_i, \psi_j\>|$ is interpreted as the mutual
coherence between $\Phi$ and $\Psi$ \cite{DonohoElad}.

For sake of concreteness, we now return to the Fourier ensemble.
Let us first set up what we are trying to prove.  Fix $\alpha > 0$,
which we shall assume to be sufficiently small.  We may take $N$ to be
large depending on $\alpha$, as the claim is vacuous when $N$ is
bounded depending on $\alpha$.  If $T$ is empty then the claim is
trivial, so from \eqref{eq:TK} we may assume that
\begin{equation}\label{tau-large}
K = \tau N \geq C \log^6 N
\end{equation}
for some (possibly) large constant $C$.

We need to prove \eqref{eq:uup}.  By self-adjointness, it would
suffice to show that with probability at least $1 -
O(N^{-\rho/\alpha})$
$$
\left|\langle F^*_{\Omega T} F_{\Omega T} f, f \rangle_{\ell_2(T)}
  - \tau \, \|f\|_{\ell_2(T)}^2\right| \leq \frac{1}{4} \tau \cdot
\|f\|_{\ell_2(T)}^2$$
for all $f \in \ell_2(T)$ and all $T$ obeying
\eqref{eq:TK}, thus $|T| \leq m$, where
\begin{equation}\label{m-def}
m := \alpha \tau N / \log^6 N.
\end{equation}

For any fixed $T$ and $f$, the above type of estimate can easily be established with high
probability by standard tools such as Lemma \ref{ms-lemma}.  The main difficulty is that
there are an exponentially large number of possible $T$ to consider, and for each fixed $T$ there
is a $|T|$-dimensional family of $f$ to consider.  The strategy is to cover the set of all $f$ of interest by
various finite nets at several scales, obtain good bounds on the size of such nets, obtain large deviation estimates
for the contribution caused by passing from one net to the net at the next scale, and sum using the union bound.

We turn to the details.  We can rewrite our goal as
$$
\left|\sum_{k \in \Omega} |\hat f(k)|^2 - \tau\,
  \|f\|_{\ell_2(T)}^2\right| \leq \frac{1}{4} \tau \cdot
\|f\|_{\ell_2(T)}^2$$
whenever $|T| \leq m$.  From Parseval's identity, this
is the same as asking that
$$
\left|\sum_{k \in \Z_N} (I_k - \tau) |\hat f(k)|^2\right| \leq \frac{1}{4}
\tau \cdot \|f\|_{\ell_2(T)}^2$$
whenever $|T| \leq m$.  Now let $U_m \subseteq \ell_2(\Z_N)$ denote the set
$$
U_m := \bigcup \{ B_{\ell_2(E)}: E \subseteq \Z_N, |E| = m \} = \{
f \in \ell_2(\Z_N): \|f\|_{\ell_2(\Z_N)} \leq 1, |\supp(f)| \leq m
\}.$$
Then the previous goal is equivalent to showing that
$$
\sup_{f \in U_m} \left|\sum_{k \in \Z_N} (I_k - \tau) |\hat f(k)|^2\right|
\leq \frac{1}{4} \tau$$
with probability at least $1 -
O(N^{-\rho/\alpha})$ for some $\rho > 0$. In fact we shall obtain the
stronger estimate
\begin{equation}\label{hatf-sum}
  \P\left( \sup_{f \in U_m} \left|\sum_{k \in \Z_N} (I_k - \tau) |\hat
      f(k)|^2\right| > \frac{1}{4} \tau\right) = O\left(\exp\left( - \frac{1}{\beta}\log^2 N\right)\right)
\end{equation}
for some constant $\beta > 0$. 

It remains to prove \eqref{hatf-sum}. The left-hand side of
\eqref{hatf-sum} is the large deviation probability of a supremum of
random sums over $U_m$.  This type of expression can be handled by
entropy estimates on $U_m$, as was done in \cite{bourgain-lambda},
\cite{bourgain-halasz}, \cite{ms}.  To follow their approach, we need
some notation.  For any $f \in \ell_2(\Z_N)$, we let $\hat f$ be its
discrete Fourier transform \eqref{eq:dft} and define the $X$ norm of
$f$ by
$$
\| f \|_X := \sqrt{N} \cdot \|\hat f \|_{\ell_\infty}.$$
Intuitively, if $f$ is a ``generic'' function bounded in
$\ell_2(\Z_N)$ we expect the $X$ norm of $f$ to be also be bounded (by
standard large deviation estimates).  We shall need this type of control
in order to apply Lemma \ref{ms-lemma} effectively.  To formalize this intuition we
shall need entropy estimates on $U_m$ in the $X$ norm. Let $B_X$ be
the unit ball of $X$ in $\ell_2(\Z_N)$.  Thus for instance $U_m$ is
contained inside the ball $\sqrt{m} \cdot B_X$, thanks to
Cauchy-Schwarz.  However we have much better entropy estimates
available on $U_m$ in the $X$ norm, which we now state.

\begin{definition}[Kolmogorov entropy]  
\label{def:entropy}
Let $X$ be a (finite-dimensional) normed vector space with norm $\|\cdot 
\|_X$, and let $B_X := \{ x \in X: \|x\|_X < 1 \}$ be the unit ball of
$X$.  If $U \subset X$ is any bounded non-empty subset of $X$ and $r >
0$, we define the \emph{covering number} $N(U, B_X, r) \in \Z^+$ to be the
least integer $N$ such that there exist elements $x_1, \ldots, x_N \in
X$ such that the balls $x_j + r B_X = \{ x \in X: \|x-x_j\|_X < r\}$,
$1 \leq j \leq N$ cover $U$, and the \emph{Kolmogorov entropy} as
\[
 {\cal E}(U, B_X, r) := \log_2(N(U, B_X, r)). 
\]
\end{definition}

\begin{proposition}\label{entropy-prop}  We have
\begin{equation}\label{entropy-universal}
{\cal E}( U_m, B_X, r ) \leq C \cdot m \log N \cdot \min(r^{-2} \log N, 1)
\end{equation}
for all $r > N^{-2}$.
\end{proposition}

This proposition is essentially contained in \cite{bourgain-lambda},
\cite{bourgain-halasz}, \cite{ms}; for sake of completeness we give a
proof of the proposition in the Appendix.  Let us assume this
proposition for the moment and conclude the proof of Lemma
\ref{teo:FourierUUP}.  Set
$$
\{J_0, \ldots, J_1\} =  \{j \in \Z : N^2 \le 2^j \le \sqrt{m}\}. 
$$
and fix $r = 2^{j}$ in Lemma \ref{entropy-prop}. By
\eqref{entropy-universal} one can find a finite subset $A_j$ of $U_m$
of cardinality
\begin{equation}\label{logo}
|A_j| \leq \exp\left(\frac{C}{1+2^{2j}} \cdot m \log^2 N\right), 
\end{equation}
such that for all $f \in U_m$, there exists $f_j \in A_j$ such that
$\|f-f_j\|_X \leq 2^j$.  Let us fix such sets $A_j$.  Then for any
$f \in U_m$, we have the telescoping decomposition
$$
f = f_{-\infty} + \sum_{J_0 \le j \le J_1} 
f_{j+1} - f_j. 
$$
where $f_j \in A_j$ and $f_{-\infty} = f - f_{J_0}$; here we have the
convention that $f_j = 0$ and $A_j = \{0\}$ if $j > J_1$. By
construction, $\|f_j - f_{j+1}\|_X \leq 2^{j+2}$, and
$\|f_{-\infty}\|_X \leq 2N^{-2}$. We write $g_j := f_{j+1} - f_j$,
thus $\|g_j\|_X\leq 2^{j+2}$.  Fix $k$ and observe the crude estimates
$$
|\hat f(k)| \le \| \hat f \|_{\ell_2}  = \|f\|_{\ell_2} = 1
\quad \text{ and } \quad |\hat f_{-\infty}(k)| \le \| f_{-\infty} \|_X
/\sqrt{N} \le 2 N^{-5/2}.
$$ 
It then follows from $|a+b|^2 \le |a|^2 +
|b|^2 + 2 |a| |b|$ that
$$
|\hat f(k)|^2 = \left|\sum_{J_0 \le j \le J_1} \hat g_j(k)\right|^2 +
O(N^{-5/2}).
$$
Multiplying this by $I_k -
\tau$ and summing, we obtain
\begin{align*}
  \left|\sum_{k \in \Z_N} (I_k - \tau) |\hat f(k)|^2 \right| & =  
  \left|\sum_{k \in \Z_N}
  (I_k - \tau)
  \left|\sum_{J_0 \le j \le J_1} \hat g_j(k)\right|^2\right| + O(N^{-3/2})\\
  &\leq 2 \sum_{J_0 \leq j \leq j' \leq J_1}
  Q(g_j,g_{j'}) + O(N^{-3/2}). 
\end{align*}
where $Q(g_j, g_{j'})$ is the nonnegative random variable
$$
Q(g_j, g_{j'}) := \left|\sum_{k \in \Z_N} (I_k - \tau) \Re(\hat
  g_j(k) \overline{\hat g_{j'}(k)})\right|.
$$
By \eqref{tau-large}, the error term $O(N^{-3/2})$ is less than $\tau/20$ (say)
if $N$ is sufficiently large.  Thus to prove \eqref{hatf-sum} it
suffices to show that
\begin{equation}\label{hatf-2}
  \P\left( \sum_{J_0 \leq j \leq j' \leq J_1} 
    \sup_{\stackrel{g_j \in A_j - A_{j+1}}{\|g_j\|_X\leq 2^{j+2}}} \sup_{\stackrel{g_{j'} \in A_{j'} - A_{j'+1}}{\|g_{j'}\|_X\leq 2^{j'+2}}}
    Q(g_j, g_{j'}) > \frac{\tau}{10}\right) = 
  O\left(\exp\left( - \frac{1}{\beta} \log^2 N\right)\right).
\end{equation}
The main difficulty is of course the presence of the suprema.  On the
other hand, the fact that the functions $g_j,g_{j'}$ are well
controlled both in entropy and in $X$ norm will allow us to handle
these suprema by relatively crude tools such as the union bound. By
the pigeonhole principle, we can bound the left-hand side of
\eqref{hatf-2} by
$$
\sum_{J_0 \leq j \leq j' \leq J_1} \P\left( \sup_{\stackrel{g_j \in
      A_j - A_{j+1}}{\|g_j\|_X\leq 2^{j+2}}} \sup_{\stackrel{g_{j'}
      \in A_{j'} - A_{j'+1}}{\|g_{j'}\|_X\leq 2^{j'+2}}} Q(g_j,
  g_{j'}) > \frac{c_0}{\log^2 N} \cdot \tau \right),
$$
for some small absolute constant $c_0$.  Since the number of pairs
$(j,j')$ is $O( \log^2 N )$, which is much smaller than
$\exp(\frac{1}{\beta}\log^2 N)$, it now suffices to show (after
adjusting the value of $\beta$ slightly) that
$$
\P\left( \sup_{\stackrel{g_j \in A_j - A_{j+1}}{\|g_j\|_X\leq
      2^{j+2}}} \sup_{\stackrel{g_{j'} \in A_{j'} -
      A_{j'+1}}{\|g_{j'}\|_X\leq 2^{j'+2}}} Q(g_j, g_{j'}) >
  \frac{c_0}{\log^2 N} \cdot \tau \right) = O\left(\exp\left( -
    \frac{1}{\beta} \log^2 N\right)\right).
$$
whenever $J_0 \leq j \leq j' \leq J_1$. 

Fix $j, j'$ as above.  From \eqref{logo} the number of possible values
of $g_{j'}$ is at most $\exp(\frac{C}{1+2^{2j'}} \cdot m \log N)$.
Thus by the union bound it suffices to show that
$$
\P\left( \sup_{\stackrel{g_j \in A_j - A_{j+1}}{\|g_j\|_X\leq
      2^{j+2}}} Q(g_j, g_{j'}) > \frac{c_0}{\log^2 N} \cdot \tau
\right) = O\left(\exp\left(-\frac{C}{1+2^{2j'}} \cdot m \log^2
    N\right)\right);
$$
for each $g_{j'} \in A_{j'} - A_{j'+1}$ with $\|g_{j'}\|_X \leq 2^{j'+2}$; note that we can absorb
the $\exp( -\frac{1}{\beta} \log^2 N)$ factor since $2^{2j'} \leq m$.

Fix $g_{j'}$.  We could also apply the union bound to eliminate the
supremum over $g_j$, but it turns out that this will lead to inferior
bounds at this stage, because of the poor $\ell_\infty$ control on
$g_{j'}$. We have to first split the frequency domain $\Z_N$ into two
sets $\Z_N = E_1 \cup E_2$, where
$$
E_1 := \{ k \in \Z_N: |\hat g_{j'}(k)| \geq \frac{C_0 2^{j} \log^2 N}{
  \sqrt{N}} \} \hbox{ and } E_2 := \{ k \in \Z_N: |\hat g_{j'}(k)| <
\frac{C_0 2^j \log^2 N}{\sqrt{N}} \}
$$
for some large absolute constant $C_0$.  Note that these sets depend
on $g_{j'}$ but not on $g_j$.  It thus suffices to show that
\begin{equation}\label{ei}
  \P\left( 
    \sup_{g_j \in A_j - A_{j+1}: \|g_j\|_X\leq 2^{j+2}} 
    Q_i(g_j, g_{j'}) > \frac{c_0}{\log^2 N} \cdot \tau \right) =  O\left(\exp\left(-\frac{C}{1+2^{2j'}} \cdot m \log^2 N\right)\right)
\end{equation}
for $i=1,2$, where we have substituted \eqref{m-def}, and
$Q_i(g_j,g_{j'})$ is the random variable
$$
Q_i(g_j, g_{j'}) := \left|\sum_{k \in E_i} (I_k - \tau) \Re(\hat
  g_j(k) \overline{\hat g_{j'}(k)})\right|.
$$
We treat the cases $i=1,2$ separately.

\textbf{Proof of \eqref{ei} when $i=1$.}  For the contribution of the
large frequencies $E_1$ we will take absolute values everywhere, which
is fairly crude but conveys the major advantage that we will be able
to easily eliminate the supremum in $g_j$. Note that since
$$ |\hat g_{j'}(k)| \leq \|g_{j'} \|_X / \sqrt N \leq 2^{j'+2} / \sqrt{N},$$
we see that this case is vacuous unless
$$ \frac{2^{j'+2}}{\sqrt{N}} \geq \frac{C_0 2^{j} \log^2 N}{\sqrt{N}},$$ or in other words, 
\begin{equation}\label{cje}
  2^{j'-j} \geq \frac{C_0 \log^2 N}{4}.
  \end{equation}
We then use the crude bound
$$ |\hat g_j(k)| \leq \| g_j \|_X / \sqrt{N} \leq 2^{j+2} / \sqrt{N}$$
and the triangle inequality to conclude
$$\sup_{g_j \in A_j - A_{j+1}} Q_1(g_j, g_{j'})
\leq \frac{2^{j+2}}{\sqrt{N}} \sum_{k \in E_1} (I_k + \tau) |\hat
g_{j'}(k)|.$$ By definition of $E_1$, we have
\begin{align*}
  \sum_{k \in E_1} 2\tau |\hat g_{j'}(k)|
  &\leq \frac{2 \tau \sqrt{N}}{C_0 2^{j} \log^2 N} \sum_{k \in \Z_N} |\hat g_{j'}(k)|^2 \\
  &= \frac{2 \tau \sqrt{N}}{C_0 2^j \log^2 N} \| g_{j'} \|_{\ell_2}^2 \\
  &\leq \frac{C \tau \sqrt{N}}{C_0 2^j \log^2 N}.
\end{align*}
Writing $I_k + \tau = (I_k - \tau) + 2\tau$, we conclude that
$$
\sup_{g_j \in A_j - A_{j+1}} Q_1(g_j, g_{j'}) \leq
\frac{2^{j+2}}{\sqrt{N}} \sum_{k \in E_1} (I_k-\tau) |\hat g_{j'}(k)|
+ \frac{C \tau}{C_0 \log^2 N}
$$
and hence to prove \eqref{ei} when
$i=1$, it would suffice (if $C_0$ is chosen sufficiently large) to show
that
$$
\P\left( \frac{2^{j+2}}{\sqrt{N}} \sum_{k \in E_1} (I_k-\tau) |\hat
  g_{j'}(k)| > \frac{c_0}{\log^2 N} \cdot \tau \right) =
O\left(\exp\left(-\frac{C}{1+2^{2j'}} \cdot m \log^2 N\right)\right).
$$
It thus suffices to show
$$
\P\left( \left|\sum_{k \in E_1} (I_k-\tau) a(k)\right| \geq \gamma \right) =
O\left(\exp\left(-\frac{C}{1+2^{2j'}} \cdot m \log^2 N\right)\right)
$$
where $a(k) := |\hat g_{j'}(k)|$ and $\gamma := \frac{c_0 \tau
  \sqrt{N}}{C 2^{j} \log^2 N}$.  Recall that
$$
\| a(k) \|_{\ell_\infty(E_1)} = O( \| g_{j'} \|_X / \sqrt{N} ) = O(
2^{j'} / \sqrt{N} )
$$
and
$$
\| a(k) \|_{\ell_2(E_1)} \leq \| \hat g_{j'} \|_{\ell_2} = \|
g_{j'} \|_{\ell_2} = O(1).
$$
We apply Lemma \ref{ms-lemma} and obtain
$$
\P\left( \left|\sum_{k \in E_1} (I_k-\tau) a(k)\right| \geq \gamma
\right ) = O\left(\exp\left(- \frac{C \gamma^2}{4 \tau + \gamma
      2^{j'}/\sqrt{N}}\right)\right).
$$
Using \eqref{cje}, we see that $\gamma 2^{j'}/\sqrt{N} \ge c
\cdot \tau$ for some absolute constant $c>0$, and conclude that
$$
\P\left( \left|\sum_{k \in E_1} (I_k-\tau) a(k)\right| \geq \gamma
\right ) = O\left( \exp\left( - \frac{C \cdot c_0 \cdot \tau
      N}{2^{j+j'} \log^2 N}\right)\right).
$$
Taking logarithms, we deduce that this contribution will be
acceptable if
$$
\frac{1}{1+2^{2j'}} \cdot m \log^2 N \le C \cdot \frac{c_0 \cdot
  \tau N}{2^{j+j'} \log^2 N}
$$
which holds (with some room to spare) thanks to \eqref{m-def}.

\textbf{Proof of \eqref{ei} when $i=2$.}  For the contribution of the small frequencies $E_2$ we use
\eqref{logo} and the union bound, and reduce to showing that
\begin{equation}\label{chn} 
\P\left( 
Q_2(g_j, g_{j'})
 > \frac{c_0}{\log^2 N} \cdot \tau \right) =    O\left(\exp\left(-\frac{C}{1+2^{2j}} m
\log^2 N)\right)\right)
\end{equation}
for any $g_j \in A_j - A_{j+1}$.  

Fix $g_j$, and set $a(k) := \Re(\hat g_j(k) \overline{\hat
  g_{j'}(k)})$; thus
$$
Q_2(g_j, g_{j'}) = \sum_{k \in E_2} (I_k - \tau) a(k).
$$
By definition of $E_2$, we have
$$
\|a(k) \|_{\ell_\infty(E_2)} \leq O( \frac{2^{j} \log^2 N}{\sqrt{N}}
\frac{\| g_j \|_X}{\sqrt{N}}) = O( \frac{2^{2j} \log^2 N}{N} )
$$
while from H\"older and Plancherel
$$
\|a(k) \|_{\ell_2(E_2)} \leq \| \hat g_{j'} \|_{\ell_2} \|
\hat g_j \|_{\ell_\infty} = \frac{\| g_{j'} \|_{\ell_2} \|
  g_j \|_X } { \sqrt{N} } = O( \frac{2^j}{\sqrt{N}} ).$$
We can apply
Lemma \ref{ms-lemma} 
to conclude
\newcommand{\cjtwo}{{c^2_{j,j'}}}
\begin{align*}
\P\left( Q_2(g_j,g_{j'}) > \frac{c_0}{\log^2 N} \cdot \tau \right) & = O( \exp(-
\frac{C \cdot \frac{c_0^2}{\log^4 N} \cdot \tau^2}{\tau [2^{2j}/N] + \tau [c_0 \log^2 N 2^{2j}/\log^2 N
  N]}))\\ & = O( \exp( - \frac{C \cdot \log^{-4} N \cdot \tau N}{2^{2j}} ) ).
\end{align*}
Taking logarithms, we thus see that this contribution will be
acceptable if
$$
\frac{1}{1+2^{2j}} \cdot m \log^2 N \le C \cdot \frac{\tau N}{2^{2j} \log^4 N}
$$
which holds thanks to \eqref{m-def}.  
This concludes the proof
of Lemma \ref{teo:FourierUUP} (assuming Proposition
\ref{entropy-prop}).  \qed

\section{`Universal' Encoding}
\label{sec:universal}

Our results interact with the agenda of coding theory.  In fact, one
can think of the process of taking random measurements as a kind of
universal coding strategy that we explain below. In a nutshell,
consider an encoder/decoder pair which would operate roughly as
follows:
\begin{itemize}
\item The Encoder and the Decoder share a collection of random vectors
  $(X_k)$ where the $X_k$'s are independent Gaussian vectors with
  standard normal entries. In practice, we can imagine that the
  encoder would send the seed of a random generator so that the
  decoder would be able to reconstruct those `pseudo-random' vectors.
  
\item {\em Encoder}. To encode a discrete signal $f$, the encoder
  simply calculates the coefficients $y_k = \<f, X_k\>$ and quantizes
  the vector $y$.
  
\item {\em Decoder}. The decoder then receives the quantized values
  and reconstructs a signal by solving the linear program
  \eqref{eq:P1}.
\end{itemize}

This encoding/decoding scheme is of course very different from those
commonly discussed in the literature of information theory. In this
scheme, the encoder would not try to know anything about the signal,
nor would exploit any special structure of the signal; it would
blindly correlate the signal with noise and quantize the
output---effectively doing very little work. In other words, the
encoder would treat each signal in exactly the same way, hence the
name ``universal encoding.'' There are several aspects of such a
strategy which seem worth exploring:

\begin{itemize}
\item {\em Robustness.} A fundamental problem with most existing
  coding strategies is their fragility vis a vis bit-loss.  Take JPEG
  2000, the current digital still-picture compression standard, for
  example. All the bits in JPEG 2000 do not have the same value and if
  important bits are missing (e.g. because of packet loss), then there
  is simply no way the information can be retrieved accurately.
  
  The situation is very different when one is using the scheme
  suggested above. Suppose for example that with a little more than
  $K$ coefficients one achieves the distortion obeying the power-law
  \begin{equation}
    \label{eq:toy}
   \| f  - f^\sharp \|^2 \lessapprox 1/K. 
\end{equation}
(This would correspond to the situation where our objects are bounded
in $\ell_1$.)  Thus receiving a little more than $K$ random
coefficients essentially allows to reconstruct a signal as precisely
as if one knew the $K$ largest coefficients.

Now suppose that in each packet of information, we have both encoded
the (quantized) value of the coefficients $y_k$ but also the label of
the corresponding coefficients $k$.  Consider now a situation in which
half of the information is lost in the sense that only half of the
coefficients are actually received. What is the accuracy of the
decoded message $f^\sharp_{50\%}$?  This essentially corresponds to
reducing the number of randomly sampled coefficients by a factor of
two, and so by \eqref{eq:toy} we see that the distortion would obey
  \begin{equation}
    \label{eq:toy2}
   \| f  - f^\sharp_{50\%} \|^2 \lessapprox 2/K  
\end{equation}
and, therefore, losses would have minimal effect.

\item {\em Security.} Suppose that someone would intercept the
  message. Then he/she would not be able to decode the message because
  he/she would not know in which random basis the coefficients are
  expressed. (In practice, in the case where one would exchange the
  seed of a random generator, one could imagine protecting it with
  standard technologies such as RSA.  Thus this scheme can be viewed as
  a variant of the standard stream cipher, based on applying a XOR
  operation between the plain text and a pseudorandom keystream, but
  with the advantage of robustness.)
    
\item {\em Cost Efficiency.} Nearly all coding scenarios work roughly
  as follows: we acquire a large number of measurements about an
  object of interest, which we then encode. This encoding process
  effectively discards most of the measured data so that only a
  fraction of the measurement is being transmitted. For concreteness,
  consider JPEG 2000, a prototype of a transform coder. We acquire a
  large number $N$ of sample values of a digital image $f$. The
  encoder then computes all the $N$ wavelet coefficients of $f$, and
  quantizes only the $B \ll N$ largest, say. Hence only a very small
  fraction of the wavelet coefficients of $f$ are actually
  transmitted.
  
  In stark contrast, our encoder makes measurements that are
  immediately used. Suppose we could design sensors which could
  actually measure the correlations $\<f, X_k\>$. Then not only the
  decoded object would be nearly as good (in the $\ell_2$-distance) as
  that obtained by knowing all the wavelet coefficients and selecting
  the largest (it is expected that the $\ell_1$-reconstruction is
  well-behaved vis a vis quantization), but we would effectively
  encode all the measured coefficients and thus, we would not discard
  any data available about $f$ (except for the quantization).
\end{itemize}

Even if one could make all of this practical, a fundamental question
remains: is this an efficient strategy? That is, for a class of
interesting signals, e.g. a class of digital images with bounded
variations, would it be possible to adapt the ideas presented in this
paper to show that this scheme does not use many more bits than what
is considered necessary? In other words, it appears interesting to
subject this compression scheme to a rigorous information theoretic
analysis. This analysis would need to address 1) how one would want to
efficiently quantize the values of the coefficients $\<f, X_k\>$ and
2) how the quantization quantitatively affects the precision of the
reconstructed signal.


\section{Discussion}
\label{sec:discussion}

\subsection{Robustness}

To be widely applicable, we need noise-aware variants of the ideas
presented in this paper which are robust against the effects of
quantization, measurement noise and modeling error, as no real-world
sensor can make perfectly accurate measurements. We view these issues
as important research topics.  For example, suppose that the
measurements $y_k = \<f, \psi_k\>$ are rounded up to the nearest
multiple of $q$, say, so that the available information is of the form
$y_k^q$ with $-q/2 \le y_k^q - y_k \le q/2$. Then we would like to
know whether the solution $f^\#$ to \eqref{eq:P1} or better, of the
variant
\[
\min_{g} \, \|g\|_{\ell_1}, \quad \text{subject to} \quad
\|F_\Omega \, g - y^q\|_{\ell_\infty} \le q/2
\]
still obeys error estimates such as those introduced in Theorem
\ref{general-lp-control}. Our analysis seems to be amenable to this
situation and work in progress shows that the quality of the
reconstruction degrades gracefully as $q$ increases. Precise
quantitative answers would help establishing the information theoretic
properties of the scheme introduced in Section \ref{sec:universal}.

\subsection{Connections with other works}
\label{connection-sec}

Our results are connected with very recent work of A.~Gilbert,
S.~Muthukrishnan, and M.~Strauss \cite{GilbertStrauss},
\cite{GilbertStraussII}.  In this work, one considers a discrete
signal of length $N$ which one would like to represent as a sparse
superposition of sinusoids. In \cite{GilbertStrauss}, the authors
develop a randomized algorithm that essentially samples the signal $f$
in the time domain $O(B^2 \text{poly}(\log N))$ times
($\text{poly}(\log N)$ denotes a polynomial term in $\log
N$) and returns a vector of approximate Fourier coefficients. They
show that under certain conditions, this vector gives, with positive
probability, an approximation to the discrete Fourier transform of
$\hat f$ which is almost as good as that obtained by keeping the
$B$-largest entries of the discrete Fourier transform of $\hat f$. In
\cite{GilbertStraussII}, the algorithm was refined so that (1) only
$O(B \text{poly}(\log N))$ samples are needed and (2) so that the
algorithm runs in $O(B \text{poly}(\log N))$ time which truly is a
remarkable feat. To achieve this gain, however, one has to sample the
signal on highly structured random grids.

Our approach is different in several aspects. First and foremost, we
are given a fixed set of nonadaptive measurements. In other words, the
way in which we stated the problem does not give us the `luxury' of
adaptively sampling the signals as in \cite{GilbertStraussII}. In this
context, it is unclear how the methodology presented in
\cite{GilbertStrauss,GilbertStraussII} would allow reconstructing the
signal $f$ from $O(B \text{poly}(\log N))$ arbitrary sampled values.
In contrast, our results guarantee that an accurate reconstruction is
possible for nearly all possible measurements sets taken from
ensembles obeying {\bf UUP} and {\bf ERP}. Second, the methodology
there essentially concerns the recovery of spiky signals from
frequency samples and do not address other setups. Yet, there
certainly is a similar flavor in the statements of their results. Of
special interest is whether some of the ideas developed by this group
of researchers might be fruitful to attack problems such as those
discussed in this article.

While finishing the write-up of this paper, we became aware of very
recent and independent work by David Donoho on a similar project
\cite{CompressedSensing}.  In that paper which appeared one month
before ours, Donoho essentially proves Theorem \ref{lp-control} for
Gaussian ensembles.  He also shows that if a measurement matrix obeys
3 conditions (CS1-CS3), then one can obtain the estimate
\eqref{approximation}. There is some overlap in methods, in particular
the estimates of Szarek \cite{Szarek2} on the condition numbers of
random matrices (CS1) also play a key role in those papers, but there
is also a greater reliance in those papers on further facts from
high-dimensional geometry, in particular in understanding the shape of
random sections of the $\ell_1$ ball (CS2-CS3). Our proofs are
completely different in style and approach, and most of our claims are
different. While \cite{CompressedSensing} only derives results for the
Gaussian ensemble, this paper establishes that other types of
ensembles such as the binary and the Fourier ensembles and even
arbitrary measurement/synthesis pairs will work as well. This is
important because this shows that concrete sensing mechanisms may be
used in concrete applications.

In a companion \cite{DecodingLP} to this paper we actually improve on
the results presented here and show that Theorem
\ref{general-lp-control} holds for general measurement ensembles
obeying the {\textbf{UUP}}. The implication for the Gaussian ensemble
is that the recovery holds with an error in \eqref{approximation} of
size at most a constant times $(K/\log(N/K))^{-r}$.



\section{Appendix: Proof of entropy estimate}
\label{entropy-estimate}

In this section we prove Proposition \ref{entropy-prop}.  
The material here is to a large extent borrowed from that in
\cite{bourgain-lambda}, \cite{bourgain-halasz}, \cite{ms}.

The entropy of the unit ball of a Hilbert space can be estimated using
the \emph{dual Sudakov inequality} of Pajor and Tomczak-Jaegerman
\cite{ptj} (See \cite{blm}, \cite{ms} for a short ``volume packing''
proof, and \cite{ms} for further discussion):

\begin{lemma}\label{entropy-lemma}\cite{ptj}  
  Let $H$ be a $n$-dimensional Hilbert space with norm $\| \cdot
  \|_{H}$, and let $B_{H}$ be the associated unit ball.  Let $e_1,
  \ldots, e_n$ be an orthonormal basis of the Hilbert space $H$, and
  let $Z_1, \ldots, Z_n \sim N(0,1)$ be i.i.d. standard Gaussian
  random variables.  Let $\| \cdot \|_Y$ be any other norm on $\C^n$.
  Then we have
  $$
   {\cal E}( B_H, B_Y, r) \leq C r^{-2} \cdot \E( \|
  \sum_{j=1}^n Z_j e_j \|_Y)^2$$
  where $C$ is an absolute constant
  (independent of $n$).
\end{lemma}

To apply this Lemma, we need to estimate the $X$ norm of certain
randomized signs.  Fortunately, this is easily accomplished:
\begin{lemma}\label{guff}  
  Let $f \in \ell_2(\ZZ_N)$ and $Z(t)$, $t \in \ZZ_N$, be i.i.d.
  standard Gaussian random variables.  Then
  $$
  \E( \| Z f \|_X ) \leq C \cdot \sqrt{\log N} \cdot \| f
  \|_{\ell_2}.$$
  The same statement holds if the $Z$'s are i.i.d.
  Bernoulli symmetric random variables ($Z(t) = \pm 1$ with equal
  probability).
\end{lemma}

\begin{proof}  Let us normalize $\|f\|_{\ell_2} = 1$.  
For any $\lambda > 0$, we have
\begin{align*}
  \P( \| Zf \|_X > \lambda ) &= \P\left( \left|\sum_{t \in \ZZ_N} Z(t) f(t)
    e^{-2\pi i t k/N}\right| > \lambda
    \hbox{ for some } k \in \ZZ_N\right) \\
  &\leq N \sup_{k \in \ZZ_N} \P\left( \left|\sum_{t \in \ZZ_N} Z(t) f(t)
    e^{-2\pi i t k/N}\right| > \lambda \right).
\end{align*}
If the $Z(t)$ are i.i.d. normalized Gaussians, then for each fixed
$k$, $\sum_{t \in \ZZ_N} Z(t) f(t) e^{-2\pi i t k/N}$ is a Gaussian
with mean zero and standard deviation $\| f\|_{\ell_2} = 1$.
Hence
$$
\P( \| Zf \|_X > \lambda ) \leq C \cdot N \cdot e^{- \lambda^2
  /2}.$$
Combining this with the trivial bound $\P( \| Zf \|_X >
\lambda ) \leq 1$ and then integrating in $\lambda$ gives the result.
The claim for i.i.d. Bernoulli variables is similar but uses
Hoeffding's inequality; we omit the standard details.
\end{proof}

Combining this lemma with Lemma \ref{entropy-lemma}, we immediately
obtain
\begin{corollary}\label{corby}  Let $E$ be a non-empty subset 
  of $\ell_2(\ZZ_N)$; note that $\ell_2(E)$ is both a Hilbert space
  (with the usual Hilbert space structure), as well as a normed vector
  space with the $X$ norm.  For all $r > 0$, we have
  $$
  {\cal E}( B_{\ell_2(E)}, B_X, r ) \leq C r^{-2} \cdot |E|
  \cdot \log N.$$
\end{corollary}

Now we turn to the set $U_m$ introduced in the preceding section.
Since the number of sets $E$ of cardinality $m$ is $\binom{N}{m} \leq
N^m$, we have the crude bound
$$
N( U_m, B_X, r ) \leq N^m \sup_{E \subseteq \ZZ_N, |E| = m}
N({\cal E}, B_X, r )
$$
and hence by Corollary \ref{corby}
\begin{equation}\label{t-bone}
{\cal E}( U_m, B_X, r ) \leq C (1+r^{-2}) \,  m \, \log N.
\end{equation}
This already establishes \eqref{entropy-universal} in the range
$C^{-1} \leq r \leq C \sqrt{\log N}$. However, this bound is quite
poor when $r$ is large.  For instance, when $r \geq m^{1/2}$ we have
\begin{equation}\label{t-zero} 
{\cal E}(U_m, B_X, r) = 0
\end{equation}
since we have $\| f \|_X < m^{1/2}$ whenever $\|f\|_{\ell_2} \le 1$
and $|\supp(f)| \leq m$.  In the regime $1 \ll r \leq m^{1/2}$ we can
use the following support reduction trick of Bourgain to obtain a
better bound:

\begin{lemma}\cite{bourgain-lambda}  
  If $r \geq C \sqrt{\log N}$ and $m \geq C$, then
  $$
  N( U_m, B_X, r ) \leq N( U_{m/2 + C \sqrt{m}}, B_X,
  \frac{r}{\sqrt{2} + C/\sqrt{m}} - C \sqrt{\log N}).$$
\end{lemma}

\begin{proof}
  Let $f \in U_m$ and $E := \supp(f)$, thus $\|f\|_{\ell_2} \le 1$ and
  $|E| \leq m$.  Let $\sigma(t) = \pm 1$ be i.i.d.  Bernoulli
  symmetric variables.  We write $f = \sigma \, f + (1-\sigma)f$.
  From Lemma \ref{guff} for Bernoulli variables we have
  $$
  \E( \| \sigma \, f \|_X ) \leq C \sqrt{\log N}$$
  and hence by Markov's
  inequality
  $$
  \P( \| \sigma \, f \|_X \geq C \sqrt{\log N} ) \leq \frac{1}{10}$$
  for a
  suitable absolute constant $C$.  Also observe that
  $$
  \| (1-\sigma) f \|_{\ell_2}^2 = \sum_{t \in E; \sigma(t) =
    -1} 4 |f(t)|^2 = 2\|f\|_{\ell_2}^2 - 2 \sum_{t \in E}
  \sigma(t) |f(t)|^2,$$ 
  and hence by Hoeffding's or Khintchine's
  inequalities and the normalization $\|f\|_{\ell_2} \le 1$
  $$
  \P( \| (1-\sigma) f \|_{\ell_2} \geq \sqrt{2} + C / \sqrt{m} )
  \leq \frac{1}{10}$$
  for a suitable absolute constant $C$.  In a
  similar spirit, we have
  $$
  \supp( (1-\sigma) f ) = |\{ t \in \supp(f): \sigma(t) = -1 \}| =
  \frac{1}{2} \supp(f) - \frac{1}{2} \sum_{t \in E} \sigma(t),$$ 
  and hence
  $$
  \P( \supp((1-\sigma) f) \geq \frac{m}{2} + C \sqrt{m} ) \leq
  \frac{1}{10}$$
  for a suitable absolute constant $C$.  Combining all
  these estimates together, we see that there exists a deterministic
  choice of signs $\sigma(t) = \pm 1$ (depending on $f$ and $E$) such that
  $$
  \| \sigma \, f\|_X \leq C \sqrt{\log N}; \,\, \| (1-\sigma) f
  \|_{\ell_2} \leq \sqrt{2} + C / \sqrt{m}; \,\, \supp((1-\sigma)f)
  \leq \frac{m}{2} + C \sqrt{m}.$$
  In particular, $f$ is within $C
  \sqrt{\log N}$ (in $X$ norm) from $(\sqrt{2} + C / \sqrt{m}) \cdot
  U_{m/2 + C \sqrt{m}}$.  We thus have
  $$
  N( U_m, B_X, r ) \leq N( (\sqrt{2} + C/\sqrt{m})
  U_{m/2 + C \sqrt{m}}, B_X, r - C \sqrt{\log N})$$
  and the claim
  follows.
\end{proof}

Iterating this lemma roughly $\log_{\sqrt{2}} \frac{r}{\sqrt{\log N}}$
times to reduce $m$ and $r$, and then applying \eqref{t-bone} once $r$
becomes comparable with $\sqrt{\log N}$, we obtain
$$
{\cal E}( U_m, B_X, r ) \leq C r^{-2} \, m \, (\log N)^2 \hbox{
  whenever } C \sqrt{\log N} \leq r \leq m^{1/2},$$
which (together
with \eqref{t-zero}) yields \eqref{entropy-universal} for all $r \geq
C \sqrt{\log N}$.

It remains to address the case of small $r$, say $N^{-2} <
r < 1/2$.  A simple covering argument (see \cite[Lemma 2.7]{ms}; the
basic point is that $B_{\ell_2(E)}$ can be covered by $O(r^{-C|E|})$
translates of $r \cdot B_{\ell_2}(E)$) gives the general inequality
$$
{\cal E}( B_{\ell_2(E)}, B_X, r ) \leq C |E| \log \frac{1}{r} +
{\cal E}( B_{\ell_2(E)}, B_X, 1 )$$
for $0 < r < 1/2$, and hence by
Corollary \ref{corby}
$$
{\cal E}( B_{\ell_2(E)}, B_X, r ) \leq C |E| \log \frac{1}{r} + C
|E| \log N .$$
Arguing as in the proof of \eqref{t-bone} we
thus have
$$
{\cal E}( U_m, B_X, r ) \leq C m (\log N + \log \frac{1}{r}),$$
which gives \eqref{entropy-universal} in the range $N^{-2} < r
< 1/2$.  This completes the proof of Proposition \ref{entropy-prop}.
\qed

\end{document}